\documentclass{amsart}
\usepackage{array}
\usepackage{multirow}
\usepackage{graphicx}
\usepackage{epsfig}
\usepackage{color,amsfonts,amsmath}
\usepackage{amssymb}
\usepackage{latexsym}
\usepackage{epsf}
\usepackage{graphicx,color,graphics}
\usepackage{amsmath,amssymb}
\usepackage{dsfont}
\usepackage[latin1]{inputenc}
\newtheorem{theorem}{Theorem}[section]

\newtheorem{remark}{Remark}

\numberwithin{equation}{section}

\begin{document}

\title[Korteweg-de Vries-Burgers equation with boundary memory]
{Boundary stabilization of the Korteweg-de Vries-Burgers equation with an infinite memory-type control and applications: A qualitative and numerical analysis}
\author[B. Chentouf, A. Guesmia, M. A. Sep\'ulveda Cort\'es and R. V\'ejar Asem]{Boumedi\`ene Chentouf$^{1*}$, Aissa Guesmia$^2$, M. A. Sep\'ulveda Cort\'es$^{3}$ and R. V\'ejar Asem$^{4}$} 
\maketitle

\pagenumbering{arabic}

\begin{center}
$^1$Kuwait University, Faculty of Science, Department of Mathematics\\
Safat 13060, Kuwait\\
$^2$Institut Elie Cartan de Lorraine, UMR 7502, Universit\'e de Lorraine\\
3 Rue Augustin Fresnel, BP 45112, 57073 Metz Cedex 03, France\\
$^3$CI$^2$MA and DIM, Universidad de Concepci\'{o}n, Concepci\'{o}n, Chile\\
$^4$Departamento de Matem\'aticas, Universidad de La Serena, La Serena, Chile
\end{center}
\begin{abstract}
This article is intended to present a qualitative and numerical analysis of well-posedness and boundary stabilization problems of the well-known Korteweg-de Vries-Burgers equation. Assuming that the boundary control is of memory type, the history approach is adopted in order to deal with the memory term. Under sufficient conditions on the physical parameters of the system and the memory kernel of the control, the system is shown to be well-posed by combining the semigroups approach of linear operators and the fixed point theory. Then, energy decay estimates are provided by applying the multiplier method. An application to the Kuramoto-Sivashinsky equation will be also given. Moreover, we present a numerical analysis based on a finite differences method and provide numerical examples illustrating our theoretical results.
\end{abstract}

\vskip0,1truecm
{\it AMS Classification.} 35B40, 35G31, 35Q35, 65M06.

\renewcommand{\thefootnote}{}
\footnotetext{E-mail addresses: $^1$boumediene.chentouf@ku.edu.kw; $^1$chenboum@hotmail.com; $^2$aissa.guesmia@univ-lorraine.fr;  $^3$masepulveda.cortes@gmail.com; $^4$rodrigo.vejar@userena.cl; $^*$Corresponding author.}

\section{Introduction}\label{sec1}

It is well-known that numerous physical phenomena exhibit both dissipation and dispersion \cite{bon1, grad, jo1, jo2, su}. This very special property is mathematically modeled by the Korteweg-de Vries-Burgers (KdVB) equation and hence has gained a considerable prominence. As a matter of fact, the KdVB and its close relatives, has been the subject of many studies (see for instance \cite{r25, bal, bon1, bo11, bo8, bub, cer20, chenmo, chengues2, den, jia1, jia, fri2, kom, li, lik, mol, po, sak1, sua, zh4}). Instead of highlighting the contribution of each of these papers, the reader is referred to \cite{chengues2} for a comprehensive discussion of this point. In turn,
it is worth mentioning that in \cite{chengues2}, the authors managed to establish well-posedness and stability outcomes for the KdVB equation with a distributed memory. In fact, it turned out that such a distributed memory term plays a role of a dissipation mechanism and hence contributes to the stability of the system. Nonetheless, boundary controls of memory-type are commonly used in practice and consequently the natural question is: when a boundary control of memory-type is applied, what is the impact on the behavior of the solutions of the KdVB equation? This is mainly the motivation of the present work. More precisely, the problem under consideration involves the third-order KdVB equation with a boundary infinite memory
\begin{equation}
\left\{
\begin{array}{ll}
\displaystyle \partial_{t} y(x,t) -\omega_0 \partial_x^2 y(x,t) +  \omega_1 \partial_x^3 y(x,t) +\omega_2 \partial_x y(x,t) + \omega_3 y(x,t) \partial_x y(x,t)=0, &  x \in  I, \, t >0,\\
y (0,t)=y(1,t)=0, & t\in \mathbb{R}_+ , \\
\partial_x y(1,t)=\omega_4 \partial_x y(0,t)+\displaystyle \int_0^{\infty}\alpha (s) \, \partial_x y(0,t-s) ds, & t\in \mathbb{R}_+ , \\
y(x,0)=y_{0} (x) , & x\in I,\\
\partial_x y(0,-t) =y_{1} (t), & t\in \mathbb{R}_+ ,
\end{array}
\right.
\label{1}
\end{equation}
where $y:I\times\mathbb{R}_+:=(0,1)\times [0,\infty) \to \mathbb{R}$ is the amplitude of the dispersive wave at the spatial variable $x$ and time $t$ \cite{ray}, $\partial_u^n$ denotes the differential operator of order $n$ with respect to $u$; that is $\dfrac{\partial^n}{\partial u^n}$. In turn, $y_0 :I\to \mathbb{R}$ and $y_1 :\mathbb{R}_+\to \mathbb{R}$ are known initial data, $\alpha:\mathbb{R}_+\to\mathbb{R_+}$ is a given function and
$\omega_i$ are real constants (physical  parameters) satisfying the following hypothesis ${\bf (H)}$:
\vskip0,1truecm
\begin{itemize}
\item The memory kernel $\alpha$ satisfies
\begin{equation}\label{f}
\alpha\in C^2 (\mathbb{R}_{+}),\quad \alpha^{\prime} < 0,\quad \alpha(0)>0,\quad \lim_{s\to\infty} \alpha(s)=0
\end{equation}
and
\begin{equation}\label{fi}
-\xi (s) \alpha^{\prime}(s)\leq \alpha^{\prime\prime} (s) \leq -\xi_0 \alpha^{\prime}(s),\quad s\in \mathbb{R}_+ ,
\end{equation}
for a positive constant $\xi_0$ and a function $\xi :\mathbb{R}_+\to \mathbb{R}_+^*:=(0,\infty)$ such that
\begin{equation}\label{xi}
\xi\in C^1 (\mathbb{R}_+),\quad \xi^{\prime}\leq 0\quad\hbox{and}\quad \int_0^{\infty} \xi (s)ds =\infty.
\end{equation}
\item The constants $\omega_0,\,\omega_1$ and $\omega_4$ satisfy
\begin{equation}
\omega_0\geq 0,\quad \omega_1 >0\quad\hbox{and}\quad\vert \omega_4\vert <1.\label{ai}
\end{equation}
\item The following relationship between $\alpha$, $\omega_1$ and $\omega_4$ holds
\begin{equation}\label{f2}
\int_0^{\infty}\frac{-\alpha^{\prime}(s)}{\xi (s)}dx:=\alpha_0 <
\frac{\omega_1\left(1-\omega_4^2\right)}{\omega_1^2 \left(1-\omega_4^2\right)+\left(1+\omega_1 \omega_4\right)^2}.
\end{equation}
\end{itemize}
\vskip0,1truecm
\begin{remark}
\begin{itemize}
\item[(i)] Typical functions $\alpha$ satisfying \eqref{f}, \eqref{fi}, \eqref{xi} and \eqref{f2} are the ones which converge exponentially to zero at infinity like \begin{equation}\label{Example1}
\alpha(s)=d_2 e^{-d_1 s},
\end{equation}
where $d_1$ and $d_2$ are positive constants satisfying ($\xi_0 =\xi (s) =d_1$)
\begin{equation}\label{Example11}
\frac{d_2}{d_1} <  \dfrac{\omega_1\left(1-\omega_4^2\right)}{\omega_1^2 \left(1-\omega_4^2\right)+\left(1+\omega_1 \omega_4\right)^2}.
\end{equation}
But \eqref{f}, \eqref{fi}, \eqref{xi} and \eqref{f2} allow $\alpha$ to have a decay rate to zero at infinity weaker than the exponential one like
\begin{equation}\label{Example2}
\alpha(s)=d_2 (1+s)^{-d_1},
\end{equation}
where $d_1 >1$ and $d_2 >0$ satisfying ($\xi_0 =d_1$ and $\xi (s) =d_1 (1+s)^{-1}$)
\begin{equation}\label{Example22}
\frac{d_2}{d_1 -1} <  \dfrac{\omega_1\left(1-\omega_4^2\right)}{\omega_1^2 \left(1-\omega_4^2\right)+\left(1+\omega_1 \omega_4\right)^2}.
\end{equation}
\vskip0,1truecm
\item[(ii)] The expression 
\begin{equation*}
\omega_4 \partial_x y(0,t)+\displaystyle \int_0^{\infty}\alpha (s) \, \partial_x y(0,t-s) ds
\end{equation*}
is viewed as a boundary control of memory-type. It is also relevant to note  that $\omega_4$ obeys the last condition in \eqref{ai} and hence one may take $\omega_4=0$. This means that, in this event, the boundary control is a purely memory-type one.
\end{itemize}
\end{remark}
\vskip0,1truecm
As mentioned earlier, our aim is to address the effect of the presence of the infinite memory term in the boundary control on the behavior of the solutions to \eqref{1}. To do so, we shall place the system in the so-called past history framework \cite{da} (see also \cite{afg} for a further discussion about the history approach and \cite{pan} for another methodology of treatment of systems with memory). The problem is shown to be well-posed as long as the hypothesis ${\bf ({H})}$ holds. Then, we prove that the memory part of the boundary control is beneficial in decaying the energy of the system under different circumstances of the kernel $\alpha$.
\vskip0,1truecm
Before providing an overview of this article, we point out that this work goes beyond the earlier one \cite{chengues2} in several respects. First, we deal with a boundary control in contrast to \cite{chengues2}, where the control is distributed. As the reader knows, it is usually more delicate (mathematically speaking) to consider a boundary control than a distributed one. Second, we mange to show that the system under consideration \eqref{1} is well-posed and its solutions are stable despite the presence of the memory term. Moreover, the desired stability property remains attainable even if the memory term is the sole action of the boundary control; that is $\omega_4=0$. The proofs are based on the multiplier method and a combination of the semigroups approach of linear operators and the fixed point theory. Third, we show that the techniques used for the KdVB system \eqref{1} can also be applied to another type of equations, namely, the Kuramoto-Sivashinsky (KS) equation. Finally, we present a numerical analysis by constructing a numerical scheme based on a finite differences method, and provide numerical examples illustrating our theoretical results.   
\vskip0,1truecm
The paper comprises six sections excluding the introduction. In Section \ref{sec2}, we put forward some preliminaries. Section \ref{sec3} is devoted to establishing the well-posedness of the KdVB system. In Section \ref{sec4}, two estimates of the energy decay rate are provided depending on the feature of the memory kernel. Indeed, it is shown that the decay of the energy corresponding of KdVB solutions is basically a transmission of the decay properties of the memory kernel. In section \ref{sec5}, we treat the KS system (\ref{11}). In Section \ref{sec6}, we give a numerical analysis for both KdVB and KS systems. Lastly, brief concluding remarks are pointed out in Section \ref{sec7}.

\section{Preliminaries}\label{sec2}

In the sequel, $\left\langle \cdot, \cdot\right\rangle$ denotes the standard real inner product in $L^2 (I)$ whose norm is $\Vert \cdot\Vert$. Then, let
$\beta=-\alpha^{\prime}$. On one hand, we deduce from \eqref{f} and \eqref{fi} that $\beta:\,\mathbb{R}_+\to \mathbb{R}_+^*$,
$\beta\in C^1 (\mathbb{R}_+)$,
\begin{equation}
-\xi_0 \beta (s)\leq \beta^{\prime} (s)\leq -\xi (s) \beta(s),\quad s\in \mathbb{R}_+\label{g01}
\end{equation}
and
\begin{equation}
\displaystyle\int_0^{\infty} \beta(s)ds=\alpha(0) >0. \label{g0}
\end{equation}
On the other hand, we infer from \eqref{g01} and \eqref{g0} that
$\beta^{\prime} < 0$, $\beta(0)>0$ as well as
\begin{equation}
\beta(0)e^{-\xi_0 s}\leq \beta(s)\leq \beta(0)e^{-\int_0^s\xi (\tau)\,d\tau}, \quad  s\in \mathbb{R}_+ ,
\label{gexp}
\end{equation}
and then, according to the last condition in \eqref{xi},
\begin{equation}
\lim_{s\to\infty}\beta (s)=0.\label{limbeta0}
\end{equation}
\vskip0,1truecm
Thereafter, following the history approach \cite{da}, we define the variable $\eta$ and its initial data $\eta_0$ as follows:
\begin{equation}\label{etaname}
\eta (t,s)= \displaystyle\int_{t-s}^t \partial_x y(0,\tau)\,d\tau\quad\hbox{and}\quad\eta_0 (s)=\displaystyle\int_{0}^s y_1 (\tau)\,d\tau,\quad t,\,s\in \mathbb{R}_+ .
\end{equation}
A formal calculation, using \eqref{1}$_5$, yields
\begin{equation}
\left\{
\begin{array}{ll}
\partial_t \eta (t,s)+\partial_s\eta (t,s) =\partial_x y(0,t) ,\quad & t,\,s>0, \vspace{0.2cm}\\
\eta (t,0) =0,\quad & t\in \mathbb{R}_+ , \vspace{0.2cm}\\
\eta (0,s) =\eta_0 (s),\quad & s\in \mathbb{R}_+.
\end{array}
\right. \label{etadef}
\end{equation}
Subsequently, in view of \eqref{f}$_4$ and \eqref{etadef}$_2$, we have
\begin{equation}  \label{etatheta}
\displaystyle\int_0^{\infty}\,\beta (s) \eta (t,s)\,ds=-\displaystyle\int_0^{\infty}\,\alpha^{\prime} (s) \eta(t,s) \,ds
=\displaystyle\int_0^{\infty}\,\alpha (s) \partial_s\eta(t,s) \,ds=\displaystyle\int_0^{\infty}\,\alpha (s) \partial_x y (0,t-s) \,ds.
\end{equation}
\vskip0,1truecm
Now, we introduce the space
\begin{equation}\label{Lg}
L_{\beta} =\left\{{\tilde\eta}:\mathbb{R}_{+}\to \mathbb{R} ,\,\displaystyle\int_0^{\infty}\,\beta (s) {\tilde\eta}^2 (s) \,ds <\infty\right\}
\end{equation}
equipped with the following inner product and corresponding norm:
\begin{equation}\label{LgInner}
\langle \eta_1 ,\eta_2\rangle_{L_\beta} =\displaystyle\int_0^{\infty}\,\beta (s) \eta_1 (s) \eta_2 (s)\,ds\quad\hbox{and}\quad\Vert{\tilde\eta}\Vert_{L_\beta}= \left(\int_{0}^{\infty}\beta(s) {\tilde\eta}^2 (s)\,ds\right)^{\frac{1}{2}} .
\end{equation}
The state space is defined by
\begin{equation}\label{Hk}
\mathcal{H} =L^2 (I)\times L_\beta
\end{equation}
and endowed with the following inner product:
\begin{equation}\label{HInner}
\langle (y_1 ,\eta_1 )^T ,(y_2 , \eta_2 )^T\rangle_{\mathcal{H}} = \langle y_1 ,y_2  \rangle +\langle \eta_1 ,\eta_2 \rangle_{L_\beta}.
\end{equation}
\vskip0,1truecm
Thereby, the problem \eqref{1} reads as
\begin{align}\label{UA}
\begin{cases}
\partial_t \Theta (\cdot, t) = \mathcal{P} \Theta (\cdot, t),\quad t>0,\\
\Theta(\cdot, 0) = \Theta_{0},
\end{cases}
\end{align}
where $\Theta =(y,\eta)^T$, $\Theta_0 =(y_0,\eta_0)^T$, the nonlinear operator $\mathcal{P}$ is defined by
\begin{equation*}
\mathcal{P} \Theta =\left(
\begin{array}{c}
\omega_0 \partial_x^2 y- \omega_1\partial_x^3 y -\omega_2 \partial_x y - \omega_3 y \partial_x y
\\
\partial_x y (0,\cdot)-\partial_s \eta
\end{array}
\right)
\end{equation*}
and its domain $\mathcal{D}(\mathcal{P})$ is given by
\begin{equation*}
\mathcal{D}(\mathcal{P})=\left\{\Theta\in \mathcal{H}; \, \mathcal{P} \Theta\in \mathcal{H} ,\,\,y\in H^1_0 (I) ,\,\, \partial_x y(1,\cdot)=\omega_4 \partial_x y(0,\cdot)+\displaystyle\int_0^{\infty}\,\beta (s) \eta (\cdot, s)\,ds,\,\,\eta (\cdot,0)=0\right\}.
\end{equation*}
\vskip0,1truecm
Finally, given $T>0$, we introduce the space
\begin{equation*}
\mathcal{M}=C \left( [0,T]; \, L^2 (I) \right) \cap L^2  \left( (0,T); \, H_0^1 (I)  \right)
\end{equation*}
whose norm will be
\begin{equation*}
\Vert \cdot \Vert_{\mathcal{M}}^2=\Vert \cdot \Vert_{C ( [0,T]; \, L^2 (I) )}^2 + \Vert \cdot \Vert_{L^2  ((0,T); \, H^1(I))}^2 .
\end{equation*}

\section{Well-posedness of the problem}\label{sec3}

The aim of this section is to prove that the problem \eqref{UA} (or equivalently the KdVB system (\ref{1})) is well-posed by means of the Fixed Point Theorem.

\subsection{The linearized system associated to (\ref{1})}\label{sub1} Note that the variables $x$, $t$ and $s$ will be omitted whenever it is unnecessary. Moreover, $C$ denotes a generic positive constant that may depend on $T$, $\alpha_0$, $\alpha (0)$, $\beta (0)$ and the parameters $\omega_i$. However,  $C$ does not depend on the initial data $\Theta_0$.
\vskip0,1truecm
The linear system associated to (\ref{1}) is
\begin{equation}
\left\{
\begin{array}{ll}
\displaystyle \partial_{t} y -\omega_0 \partial_x^2 y + \omega_1 \partial_x^3 y+\omega_2 \partial_x y  =0, & x \in I,\, t>0,\\
\partial_t \eta (t,s)+\partial_s \eta (t,s) -\partial_x y(0,t)=0,\quad & t,\,s>0 ,\\
\eta (t,0) =0,\quad & t\in \mathbb{R}_+ ,\\
y(0,t) =y(1,t) =0, & t\in \mathbb{R}_+,\\
\partial_x y(1,t)=\omega_4 \partial_x y(0,t)+\displaystyle\int_0^{\infty}\,\beta (s) \eta (t, s)\,ds, & t\in \mathbb{R}_+,\\
y(x,0) =y_0 (x), & x\in I,\\
\partial_x y(0,-t) =y_1 (t), & t\in \mathbb{R}_+,\\
\eta (0, s)=\eta_0 (s)=\displaystyle\int_{0}^s y_1 (\tau)\,d\tau, & s \in \mathbb{R}_+.
\end{array}
\right.
\label{lin}
\end{equation}
Taking $\Theta_0=(y_0,\eta_0)^T$, the latter takes the abstract form in
$\mathcal{H}$
\begin{equation}
\left\{
\begin{array}{ll}
\partial_t \Theta (\cdot, t)=\mathcal{A} \Theta (\cdot, t),\quad t>0,\\
\Theta (\cdot, 0)=\Theta_{0},
\end{array}
\right. \label{si}
\end{equation}
in which $\mathcal{A}$ is the linear operator defined by
\begin{equation}
\mathcal{A} \Theta =
\begin{array}{l}
 \left(
\begin{array}{c}
\omega_0 \partial_x^2 y- \omega_1\partial_x^3 y-\omega_2 \partial_x y \\
\\
\partial_x y(0,\cdot)-\partial_s \eta
\end{array}
\right)
\end{array} \label{1.62}
\end{equation}
and its domain ${\mathcal{D}}\left(\mathcal{A}\right)$ is given by
\begin{equation*}
{\mathcal{D}}\left(\mathcal{A}\right)=\left\{
\Theta \in \mathcal{H}; \, \mathcal{A} \Theta \in \mathcal{H} ,\,\,y\in H^1_0 (I) ,\,\, \partial_x y(1,\cdot)=\omega_4 \partial_x y(0,\cdot)+\displaystyle\int_0^{\infty}\,\beta (s) \eta (\cdot, s)\,ds,\,\,\eta (\cdot,0)=0\right\}.
\end{equation*}
\vskip0,1truecm
\begin{theorem}\label{th310}
Assume that ${\bf(H)}$ holds. Then, we have:
\vskip0,1truecm
\begin{itemize}
\item[(i)] The linear operator $\mathcal{A}$ defined by (\ref{1.62}) generates a $C_{0}$-semigroup of contractions $e^{t\mathcal{A}}$ on $\mathcal{H}$. Therefore, for any initial data $\Theta_{0}=(y_0 ,\eta_0 )^T\in{\mathcal{D}}(\mathcal{A})$, the problem (\ref{si}) has a unique classical solution $\Theta =(y ,\eta )^T$ satisfying
\begin{equation}
\Theta\in C(\mathbb{R}_+;{\mathcal{D}}(\mathcal{A})) \cap C^{1}(\mathbb{R}_+;\mathcal{H}). \label{regularity}
\end{equation}
Notwithstanding, if $\Theta_{0}\in \mathcal{H}$, then (\ref{si}) admits a mild solution \begin{equation}
\Theta \in  C(\mathbb{R}_+;\mathcal{H}). \label{regularity0}
\end{equation}
\vskip0,1truecm
\item[(ii)]  For any $\Theta_{0}\in \mathcal{H}$ and $T>0$, we have the estimates
\begin{equation}
\Vert \partial_x y(0,\cdot) \Vert_{L^2 (0,T)}^2 \leq C\Vert \Theta_0\Vert_{\mathcal{H}}^2 \quad\hbox{and}\quad
\Vert \partial_x y \Vert_{L^2 ((0,T);L^2 (I) )}^2  \leq  C  \Vert  \Theta_0\Vert_{\mathcal{H}}^2 , \label{6}
\end{equation}
for some constant $C>0$. Finally, the map
\begin{equation}
\Xi: \Theta_0 \in \mathcal{H} \mapsto \Xi (\Theta_0)=\Theta (\cdot)=e^{\cdot\mathcal{A}}\Theta_0\in \mathcal{M} \times C \left( [0,T]; \, L_\beta \right)
\end{equation}
is continuous.
\end{itemize}
\label{th1}
\end{theorem}
\vskip0,1truecm
\begin{proof}
(i) Given $\Theta\in{\mathcal{D}}({\mathcal{A}})$, we infer from \eqref{LgInner}, \eqref{HInner}, \eqref{1.62} and simple integration by parts that
\begin{eqnarray*}
\langle{\mathcal{A}}\Theta,\Theta \rangle_{\mathcal{H}} &=& \int_0^{1} y\left(\omega_0 \partial_x^2 y- \omega_1\partial_x^3 y-\omega_2 \partial_x y\right)dx +\int_0^{\infty}\beta (s)\left(\partial_x y(0,\cdot)-\partial_s \eta\right)\eta ds\\
&=& \omega_0 \left[ y \partial_x y \right]_0^1 -\omega_0 \Vert \partial_x y \Vert^2 - \omega_1\left[ y \partial_x^{2}y \right]_0^1 +\frac{\omega_1}{2} \left[ (\partial_x y)^2\right]_0^1 -\frac{\omega_2}{2} \left[y^2\right]_0^1
\\
&& + \partial_x y(0,t)\int_0^{\infty} \beta(s)\eta ds -\frac{1}{2}
\left[ \beta (s)\eta^2\right]_0^{\infty} +\frac{1}{2}\int_0^{\infty}\beta^{\prime}(s)\eta^2 ds ,
\end{eqnarray*}
where the last integral is well defined because $\vert\beta^{\prime}\vert=-\beta^{\prime}\leq \xi_0\beta$ and $\eta\in L_{\beta}$. Then, using \eqref{limbeta0} and the boundary conditions in \eqref{lin}, we observe that
\begin{eqnarray}\label{di1}
\langle{\mathcal{A}}\Theta,\Theta \rangle_{\mathcal{H}} &=& -\omega_0\Vert \partial_x y \Vert^2 -\frac{\omega_1}{2} \left(1-\omega_4^2\right)(\partial_{x}y)^2 (0,\cdot) +\frac{\omega_1}{2}\left(\int_0^{\infty} \beta(s) \eta ds\right)^2 \\
&& + (1+\omega_1\omega_4)\partial_x y(0,\cdot)\int_0^{\infty} \beta(s)\eta ds +\frac{1}{2}\int_0^{\infty}\beta^{\prime}(s)\eta^2 ds. \nonumber
\end{eqnarray}
Using H\"older's inequality, \eqref{f2}, \eqref{g01} and the right inequality in \eqref{g01}, we see that
\begin{equation} \label{di2}
\left(\int_0^{\infty} \beta(s) \eta ds\right)^2\leq \left(\int_0^{\infty} \frac{\beta(s)}{\xi (s)} ds\right)\int_0^{\infty} \xi (s)\beta(s) \eta^2 ds \leq -\alpha_0\int_0^{\infty} \beta^{\prime}(s) \eta^2 ds.
\end{equation}
On the other hand, using Young's inequality, we find that
\begin{equation} \label{di3}
\left\vert\partial_x y(0,\cdot)\int_0^{\infty} \beta(s)\eta ds\right\vert\leq \frac{1}{2}\left[\epsilon (\partial_x y)^2 (0,\cdot) +\frac{1}{\epsilon}\left(\int_0^{\infty} \beta(s) \eta ds\right)^2\right],
\end{equation}
for any $\epsilon>0$. Combining \eqref{di1}, \eqref{di2} and \eqref{di3}, we obtain
\begin{eqnarray}
\langle{\mathcal{A}}\Theta ,\Theta \rangle_{\mathcal{H}}
&\le& -\omega_0\Vert \partial_x y \Vert^2 -\kappa (\partial_x y)^2 (0,\cdot) +\kappa \displaystyle \int_0^{\infty} \beta^{\prime} (s)\eta^2 \,ds, \label{di}
\end{eqnarray}
where
\begin{equation}\label{mu0}
\kappa =\frac{1}{2} \min\left\{\omega_1\left(1- \omega_4^2\right) -\epsilon \vert 1+ \omega_1\omega_4 \vert, 1-\alpha_0\left(\omega_1+\frac{1}{\epsilon}\vert 1+ \omega_1\omega_4 \vert\right)\right\}.
\end{equation}
It is noteworthy that if $\omega_1\omega_4 =-1$, then we infer from \eqref{f2}  that $\alpha_0 <\frac{1}{\omega_1}$. Hence, thanks to \eqref{ai}, we have
$\kappa>0$. In turn, if $\omega_1\omega_4 \ne -1$, then we still have $\alpha_0 <\frac{1}{\omega_1}$ by virtue of \eqref{f2}. Thereby, we choose $\epsilon$ as follows:
\begin{equation*}
\frac{\alpha_0\vert 1+ \omega_1\omega_4\vert}{1 - \omega_1\alpha_0} < \epsilon <\frac{\omega_1 \left(1- \omega_4^2\right)}{\vert 1+ \omega_1\omega_4\vert},
\end{equation*}
which is possible in view of \eqref{f2}. Thus, also in this case, we get that $\kappa>0$. Besides, owing to \eqref{ai} and \eqref{g01}, the linear operator $\mathcal{A}$ is dissipative.
\vskip0,1truecm
Moreover, one can readily verify that the adjoint operator of $\mathcal{A}$ is defined by: for $\Psi =\left(y^* ,\eta^*\right)^T$,
\begin{equation}
\mathcal{A}^* \Psi=
\begin{array}{l}
 \left(
\begin{array}{c}
\omega_0\partial_x^2 y^* +\omega_1\partial_x^3 y^* +\omega_2\partial_x y^* \\
\\
\omega_1\partial_x y^* (1,\cdot)+ \partial_s \eta^* + \frac{\beta^{\prime}}{\beta}\eta^*
\end{array}
\right)
\end{array}
\label{1.62n}
\end{equation}
with domain
\begin{equation*}
{\mathcal{D}}\left(\mathcal{A}^*\right)=\left\{
\Psi \in \mathcal{H}; \, \mathcal{A}^* \Psi \in \mathcal{H},\, y^*\in H^1_0 (I) ,\, \partial_x y^* (0,\cdot)=\omega_4 \partial_x y^* (1,\cdot)+\frac{1}{\omega_1}\displaystyle \int_0^{\infty} \beta(s)\eta^* (\cdot,s) \,ds,\, \eta^* (\cdot,0)=0\right\}.
\end{equation*}
Indeed, direct integrations by parts lead to
\begin{eqnarray*}
\langle{\mathcal{A}}\Theta,\Psi \rangle_{\mathcal{H}} &=& \int_0^{1} y^*\left(\omega_0 \partial_x^2 y- \omega_1\partial_x^3 y-\omega_2 \partial_x y\right)dx +\int_0^{\infty}\beta (s)\left(\partial_x y(0,\cdot)-\partial_s \eta\right)\eta^* ds\\
&=& \omega_0 \left[ y^* \partial_x y -y \partial_x y^* \right]_0^1 +\omega_0
\int_0^{1} y\partial_x^2 y^* dx- \omega_1\left[ y^* \partial_x^{2} y -\partial_x y^* \partial_x y +y \partial_x^{2} y^*\right]_0^1 \\
&& +\omega_1 \int_0^{1} y\partial_x^3 y^* dx -\omega_2 \left[yy^*\right]_0^1 + \omega_2
\int_0^{1} y\partial_x y^* dx \\
&& + \partial_x y(0,t)\int_0^{\infty} \beta(s)\eta^* ds -
\left[ \beta (s)\eta\eta^*\right]_0^{\infty} +\int_0^{\infty}\beta(s)\eta \left(\partial_s\eta^* +\frac{\beta^{\prime} (s)}{\beta (s)}\eta^*\right)ds .
\end{eqnarray*}
Thereafter, exploiting \eqref{limbeta0} and the Dirichlet boundary conditions in
${\mathcal{D}}\left(\mathcal{A}\right)$ and ${\mathcal{D}}\left(\mathcal{A}^*\right)$, we find that
\begin{eqnarray*}
\langle{\mathcal{A}}\Theta,\Psi \rangle_{\mathcal{H}} &=& \int_0^{1} y \left(\omega_0 \partial_x^2 y^* +\omega_1\partial_x^3 y^* +\omega_2 \partial_x y^*\right)dx +\int_0^{\infty}\beta (s)\eta\left(\partial_s \eta^* +\frac{\beta^{\prime} (s)}{\beta (s)}\eta^*\right) ds\\
&& +\partial_x y(0,t)\int_0^{\infty} \beta(s)\eta^* ds+ \omega_1\left[\partial_x y (1,\cdot)\partial_x y^* (1,\cdot)-\partial_x y (0,\cdot)\partial_x y^* (0,\cdot)\right].
\end{eqnarray*}
Now, making use of the Neumann boundary condition in
${\mathcal{D}}\left(\mathcal{A}\right)$, we arrive at
\begin{equation*}
\langle{\mathcal{A}}\Theta,\Psi \rangle_{\mathcal{H}} = \langle{\Theta,\mathcal{A}}^*\Psi \rangle_{\mathcal{H}}+\partial_x y (0,\cdot) \left[\int_0^{\infty}\beta (s)\eta^* ds-\omega_1\partial_x y^* (0,\cdot)+\omega_1\omega_4\partial_x y^* (1,\cdot)\right],
\end{equation*}
which, together with the Neumann boundary condition in ${\mathcal{D}}\left(\mathcal{A}^*\right)$, implies that
\begin{equation*}
\langle{\mathcal{A}}\Theta,\Psi \rangle_{\mathcal{H}} = \langle{\Theta,\mathcal{A}}^*\Psi \rangle_{\mathcal{H}}.
\end{equation*}
Whereupon, the definitions of ${\mathcal{A}}^*$ and its domain ${\mathcal{D}}\left(\mathcal{A}^*\right)$ are justified.
\vskip0,1truecm
Subsequently, one can show analogously to \eqref{di} that
\begin{eqnarray*}
\left\langle{\mathcal{A}^*} \Psi,\Psi \right\rangle_{\mathcal{H}} &\le & -\omega_0\Vert \partial_x y^* \Vert^2 -\kappa^*(\partial_x y^* )^2 (1,\cdot) +\kappa^* \displaystyle \int_0^{\infty} \beta^{\prime} (s)\eta^{*2} \,ds,
\end{eqnarray*}
for any $\Psi =(y^* ,\eta^*) \in  {\mathcal{D}}\left(\mathcal{A}^*\right)$, where
\begin{equation*}
\kappa^* =\frac{1}{2} \min\left\{\omega_1\left(1- \omega_4^2\right) -\epsilon \vert \omega_1 +\omega_4 \vert, 1-\alpha_0\left(\frac{1}{\omega_1}+\frac{1}{\epsilon}\vert \omega_1 +\omega_4 \vert\right)\right\}
\end{equation*}
and $\epsilon =1$ if $\omega_1 +\omega_4 =0$. However, if $\omega_1 +\omega_4 \ne 0$, then $\epsilon >0$ is chosen such that
\begin{equation*}
\frac{\alpha_0 \omega_1\vert\omega_1 +\omega_4\vert}{\omega_1 -\alpha_0} < \epsilon <\frac{\omega_1 \left(1- \omega_4^2\right)}{\vert \omega_1 +\omega_4\vert}.
\end{equation*}
Note that we infer from \eqref{f2} that $\alpha_0 <\omega_1$ and hence $\epsilon$ is well-defined. Thus, $\kappa^* >0$ and $\mathcal{A}^*$ is also dissipative.
\vskip0,1truecm
Lastly, since $\mathcal{A}$ is  a closed and densely defined operator, the first part (i) of Theorem \ref{th1} follows from semigroups theory of linear operators \cite{br, Pa}.
\vskip0,1truecm
(ii) Picking up $\Theta_{0} \in \mathcal{H}$ and using the contraction of the semigroup $e^{t\mathcal{A}}$, we obtain
\begin{equation}
\Vert \Theta (\cdot, t)\Vert_{\mathcal{H}}^2 =\Vert  y(\cdot, t) \Vert^2 +\Vert \eta (t,\cdot)\Vert_{L_\beta}^2  \leq  \Vert y_0 \Vert^2 + \Vert \eta_0 \Vert_{L_\beta}^2 =\Vert \Theta_0\Vert_{\mathcal{H}}^2, \quad t \in [0,T].
\label{con}
\end{equation}
\vskip0,1truecm
Next, let $p:I\times\mathbb{R}_+\to\mathbb{R}$ and $q:\mathbb{R}_+\times\mathbb{R}_+\to\mathbb{R}$ be two smooth functions. Consider $\Theta_0\in {\mathcal{D}}(\mathcal{A})$ and the solution $\Theta$ of \eqref{si} with the regularity \eqref{regularity} (a standard argument of density allows to extend the next results to solutions stemmed from $\Theta_0 \in \mathcal{H}$). Then, multiplying (\ref{lin})$_1$ by $py$, integrating by parts over
$[0,T]\times I$ and using the boundary conditions in \eqref{lin}, we obtain
\begin{eqnarray}
&& \displaystyle \int_0^1 \left( p(x,T) y^2(x,T) -p(x,0) y_0^2 (x) \right) dx - \int_0^T  \int_0^1 \left( \partial_t p+\omega_0 \partial_x^2 p +\omega_1\partial_x^3 p +\omega_2 \partial_x p   \right) y^2 (x,t)~dx dt \nonumber\\
&& + \int_0^T  \int_0^1 \left( 2\omega_0 p +3\omega_1 \partial_x p \right) \left( \partial_x y \right)^2 (x,t) dx dt +\omega_1\int_0^T   \left( p(0,t) - \omega_4^2 p(1,t)  \right) \left( \partial_x y \right)^2 (0,t) dt\nonumber\\
&& =\omega_1\int_0^{T}  p(1,t) \left( \int_0^{\infty} \beta(s)  \eta ds \right)^2 dt+ 2 \omega_1\omega_4 \int_0^{T}  p(1,t) \partial_x y(0,t)  \int_0^{\infty} \beta(s)  \eta ~ ds dt.\label{ux}
\end{eqnarray}
Then, we take the inner product in $L_\beta$ of \eqref{lin}$_2$ with
$q \eta$ and then integrate over $[0,T]$ to get
\begin{eqnarray}
&& \displaystyle \int_0^{\infty} \beta(s) \left(q(T,s) \eta^2 (T,s)- q(0,s) \eta_0^2(s)\right) ds -  \int_0^T \int_0^{\infty} \beta(s) \eta^2 \partial_t q ds dt  \nonumber\\
&& \quad =2 \int_0^T \partial_x y(0,t) \int_0^{\infty} \beta(s) q\eta dsdt+
\int_0^T \int_0^{\infty} \beta^{\prime}(s)  q\eta^2 ds dt+\int_0^T \int_0^{\infty} \beta(s)\eta^2 \partial_s q  ds dt.
\label{etat}
\end{eqnarray}
Choosing $p\equiv 1$ and $q\equiv1$ in \eqref{ux} and \eqref{etat}, respectively, and adding the obtained formulas, we have
\begin{eqnarray}
&& 2 \omega_0 \Vert  \partial_x y \Vert_{L^2 ((0,T); L^2 (I))}^2 + \omega_1 (1-\omega_4^2) \Vert  \partial_x y(0,\cdot)\Vert_{L^2 (0,T)}^2 = \Vert  y_0 \Vert^2 - \Vert y(\cdot,T) \Vert^2 + \Vert  \eta_0 \Vert_{L_\beta}^2 -\Vert  \eta (T,\cdot) \Vert_{L_\beta}^2  \nonumber\\
&& \;\; +\int_{0}^{T} \int_0^{\infty} \beta^{\prime} (s) \eta^2 ds dt+\omega_1
\int_{0}^{T}\left(\int_0^{\infty} \beta(s) \eta ds\right)^2 dt+2 (1+\omega_1\omega_4) \int_{0}^{T} \partial_x y(0,t) \int_0^{\infty} \beta(s)\eta ds dt.
\label{simple}
\end{eqnarray}
Inserting \eqref{di2} and \eqref{di3} in  \eqref{simple} yields
\begin{equation}\label{es1}
2\omega_0\Vert  \partial_x y \Vert_{L^2 ((0,T); L^2 (I))}^2 + 2\kappa \Vert  \partial_x y(0,\cdot)\Vert_{L^2 (0,T)}^2 \le \Vert  y_0 \Vert^2 + \Vert  \eta_0 \Vert_{L_\beta}^2 =\Vert  \Theta_0\Vert_{\mathcal{H}}^2 ,
\end{equation}
where $\kappa$ is the positive constant given by \eqref{mu0}.
\vskip0,1truecm
It is clear, from \eqref{con} and \eqref{es1}, that if $\omega_0 >0$, then \eqref{6} holds and the map $\Xi$ is continuous. In turn, if $\omega_0 =0$, then we have only the first estimate of \eqref{6}. In order to get the second one, let $p(x,t)=x$ in \eqref{ux}, which gives
\begin{eqnarray}
&&\displaystyle \int_0^1 x \left( y^2 (x,T)- y_0^2 (x) \right) dx -\omega_2 \int_0^T  \Vert  y \Vert^2 dt+ 3 \omega_1\int_0^T  \Vert  \partial_x y \Vert^2 dt - \omega_1\omega_4^2 \Vert  \partial_x y(0,\cdot)\Vert_{L^2 (0,T)}^2 \nonumber\\
&& =\omega_1\int_{0}^{T}\left(\int_0^{\infty} \beta(s) \eta ds\right)^2 dt +2 \omega_1\omega_4 \int_{0}^{T} \partial_x y(0,t) \int_0^{\infty} \beta(s)\eta ds dt.
\label{px}
\end{eqnarray}
Amalgamating \eqref{di2}, \eqref{di3} and \eqref{px} and using the left inequality in \eqref{g01}, we get
\begin{eqnarray}
\displaystyle 3 \omega_1\int_0^T  \Vert \partial_x y \Vert^2 dt &\leq&  \Vert y_0\Vert^2 + \omega_2 \int_0^T  \Vert y \Vert^2 dt +\left(\vert\omega_1\omega_4\vert +\omega_1\omega_4^2\right) \Vert\partial_x y(0,\cdot)\Vert_{L^2 (0,T)}^2 \nonumber\\
&&+\alpha_0\left( \omega_1+\vert\omega_1\omega_4\vert\right)\xi_0 \int_0^{T} \int_0^{\infty} \beta (s) \eta^2 ds dt. \label{k0}
\end{eqnarray}
Lastly, it suffices to combine the first estimate of \eqref{6} and \eqref{con} with \eqref{k0} to obtain
\begin{equation}
\displaystyle\int_0^T  \Vert \partial_x y \Vert^2 dt \leq C\left( \Vert y_0 \Vert^2 + \Vert \eta_0 \Vert_{L_\beta}^2 \right)=C\Vert \Theta_0 \Vert_{\mathcal{H}}^2, \label{k0f}
\end{equation}
for some positive constant $C$. Hence the second estimate of \eqref{6} is derived and the continuity of $\Xi$ follows from \eqref{con}.
\end{proof}
\bigskip
Subsequently, let us define the energy $E$ of \eqref{1} (and also \eqref{lin}) by
\begin{equation}\label{E}
E (t) =\frac{1}{2}\Vert \Theta(\cdot ,t)\Vert_{\mathcal{H}}^2 =\frac{1}{2} \left(\Vert y(\cdot ,t)\Vert^2 +\displaystyle\int_0^{\infty} \beta(s)\eta^2 (t,s)\, ds\right),\quad t\in \mathbb{R}_+ .
\end{equation}
Then, multiplying \eqref{1}$_1$ by $y$ and integrating over $I$, we get
\begin{equation*}
\partial_t \left(\frac{1}{2}\Vert y\Vert^2\right) -\omega_0 \left[ y \partial_x y \right]_0^1 +\omega_0 \Vert \partial_x y \Vert^2 + \omega_1\left[ y \partial_x^{2}y \right]_0^1 -\frac{\omega_1}{2} \left[ (\partial_x y)^2\right]_0^1 +\frac{\omega_2}{2} \left[y^2\right]_0^1 +\frac{\omega_3}{3} \left[y^3\right]_0^1 =0.
\end{equation*}
By virtue of the boundary conditions \eqref{1}$_2$, \eqref{1}$_3$ and \eqref{etatheta}, the latter becomes
\begin{equation}
\partial_t \left(\frac{1}{2}\Vert y\Vert^2\right) =-\omega_0\Vert \partial_x y \Vert^2 -\frac{\omega_1}{2} (\partial_{x}y)^2 (0,t) +\frac{\omega_1}{2}\left(\omega_4 \partial_x y(0,t)+\int_0^{\infty} \beta(s) \eta ds\right)^2 . \label{Equa1}
\end{equation}
Multiplying \eqref{etadef}$_1$ by $\beta(s)\eta$ and integrating on
$\mathbb{R}_+$, we get
\begin{equation*}
\partial_t \left(\frac{1}{2}\int_0^{\infty}\beta(s) \eta^2 ds\right) +\frac{1}{2}\int_0^{\infty}\beta(s) \partial_s \left(\eta^2 \right) ds = \partial_x y(0,t)\int_0^{\infty} \beta(s) \eta ds.
\end{equation*}
Thanks to an integration by parts and using \eqref{limbeta0} and \eqref{etadef}$_2$, we arrive at
\begin{equation}\label{Equa2}
\partial_t \left(\frac{1}{2}\int_0^{\infty}\beta(s)\eta^2 ds\right)=\frac{1}{2}\int_0^{\infty}\beta^{\prime}(s)\eta^2 ds + \partial_x y(0,t)\int_0^{\infty} \beta(s)\eta ds.
\end{equation}
Combining \eqref{Equa1} and \eqref{Equa2}, we have
\begin{eqnarray}\label{EprimeFirst}
E^{\prime} (t) &=& -\omega_0 \Vert \partial_x y \Vert^2 -\frac{\omega_1}{2} \left(1-\omega_4^2\right)(\partial_{x}y)^2 (0,t) +\frac{\omega_1}{2}\left(\int_0^{\infty} \beta(s) \eta ds\right)^2 \\
&& +\frac{1}{2}\int_0^{\infty}\beta^{\prime}(s)\eta^2 ds + (1+\omega_1\omega_4 )\partial_x y(0,t)\int_0^{\infty} \beta(s)\eta ds.\nonumber
\end{eqnarray}
In light of \eqref{di2} and \eqref{di3}, we see that
\begin{equation}
E^{\prime} (t) \leq -\omega_0 \Vert \partial_x y \Vert^2 -\kappa (\partial_{x}y)^2 (0,t) +\kappa\int_0^{\infty} \beta^{\prime}(s)\eta^2 ds,\label{Eprime0}
\end{equation}
where $\kappa$ is given by \eqref{mu0}. Lastly, by virtue of \eqref{ai}, \eqref{g01} and \eqref{Eprime0}, we can claim that the energy $E$ is non-increasing along the solutions of the system \eqref{1} and also \eqref{lin}.

\subsection{A non-homogeneous linear system associated to (\ref{1})}\label{sub2} Consider now the linear system (\ref{lin}) with an additional source term $\sigma (x,t)$
\begin{equation}
\left \{
\begin{array}{ll}
\displaystyle \partial_{t} y -\omega_0 \partial_x^2 y + \omega_1\partial_x^3 y+\omega_2 \partial_x y =\sigma(x,t), & x \in I,\, t>0,\\
\partial_t \eta (t,s)+\partial_s \eta (t,s)-\partial_x y(0,t)=0,\quad & t,\,s>0 ,\\
\eta (t,0) =0,\quad & t\in \mathbb{R}_{+} ,\\
y(0,t) =y(1,t) =0, & t \in \mathbb{R}_{+} ,\\
\partial_x y(1,t)=\omega_4 \partial_x y(0,t) +\displaystyle\int_0^{\infty} \beta (s)\eta^2 (t,s) ds, & t\in \mathbb{R}_{+} ,\\
y(x,0) =y_0 (x), & x \in I ,\\
\partial_x y(0,-t) =y_1 (t), & t\in \mathbb{R}_{+} ,\\
\end{array}
\right.
\label{fxt}
\end{equation}
with some initial data $\Theta_{0}$. We have the following result:
\vskip0,1truecm
\begin{theorem} Assume that  ${\bf(H)}$ holds. Given $T>0$, we have
\vskip0,1truecm
\begin{itemize}
\item[(i)] If $\Theta_{0}=(y_0 ,\eta_0 )^T\in \mathcal{H}$ and $\sigma\in L^1 ((0,T);L^2 (I))$, then there exists a unique mild solution $\Theta=(y,\eta )^T$ of (\ref{fxt}) such that
$\Theta\in \mathcal{M} \times C([0,T]; L_\beta )$,
\begin{equation}
\Vert \Theta \Vert_{C([0,T]; \mathcal{H})}^2  \leq  C_0 \left(\Vert \Theta_{0}\Vert_{\mathcal{H}}^2 +\Vert \sigma\Vert_{L^1 ((0,T);L^2 (I))}^2\right) \label{80}
\end{equation}
and
\begin{equation}
\Vert y \Vert_{\mathcal{M}}^2  \leq  C_1 \left(  \Vert \Theta_{0} \Vert_{\mathcal{H}}^2 + \Vert \sigma\Vert_{L^1 ((0,T);L^2 (I))}^2 \right),\label{81}
\end{equation}
for some constants $C_0 ,\,C_1 >0$ independent of $\Theta_{0}$ and
$\sigma$.
\vskip0,1truecm
\item[(ii)]
 Given ${\tilde y} \in L^2 ((0,T);H^1 (I))$, we have ${\tilde y}  \partial_x {\tilde y}  \in L^1 ((0,T);L^2 (I))$ and the map
\begin{equation*}
\Lambda: {\tilde y} \in L^2 ((0,T);H^1 (I)) \mapsto {\tilde y} \partial_x {\tilde y} \in L^1 ((0,T);L^2 (I))
\end{equation*}
is continuous.
\end{itemize}
\label{th2}
\end{theorem}
\vskip0,1truecm
\begin{proof}
(i) Thanks to the contraction of the semigroup $e^{t\mathcal{A}}$ and the fact that $\sigma \in L^1 ((0,T);L^2 (I) )$, the existence and uniqueness results follow \cite{Pa}. Next, it suffices to show the statement (i) for an initial data $\Theta_{0}$ in $\mathcal{D} (\mathcal{A})$ as previously done. Furthermore, let us define energy of (\ref{fxt}) by \eqref{E}. Next, arguing as for \eqref{Eprime0}, we get
\begin{equation}
E^{\prime}(t) \le -\omega_0 \Vert \partial_x y \Vert^2 -\kappa (\partial_{x}y)^2 (0,\cdot) +\kappa \int_0^{\infty} \beta^{\prime}(s) \eta^2 ds + \int_{0}^{1} y(x,t)  \sigma(x,t) dx \le \Vert y\Vert\Vert\sigma\Vert , \label{fxt1}
\end{equation}
where we used Cauchy-Schwarz inequality. Now, we integrate \eqref{fxt1} to find
\begin{equation*}
\Vert\Theta (\cdot,t)\Vert_{\mathcal{H}}^2 \le \Vert\Theta_0\Vert_{\mathcal{H}}^2 + 2\int_{0}^{t}\Vert y(\cdot,\tau )\Vert\Vert\sigma (\cdot,\tau )\Vert d\tau \le \Vert\Theta_0\Vert_{\mathcal{H}}^2 +2\max_{\tau\in [0,T]}\Vert y(\cdot,\tau)\Vert \int_{0}^{T}\Vert\sigma (\cdot,\tau )\Vert d\tau,
\end{equation*}
therefore, using Young's inequality to get
\begin{equation}
\displaystyle \max_{t \in [0,T]} \Vert \Theta (\cdot, t) \Vert_{\mathcal{H}}^2  \leq  \Vert \Theta_0\Vert_{\mathcal{H}}^2  + \dfrac{1}{2} \displaystyle \max_{t \in [0,T]} \Vert y(\cdot,t)\Vert^2  +2\Vert \sigma\Vert_{L^1 ((0,T);L^2 (I) )}^2. \label{83}
\end{equation}
Thus, the estimate (\ref{80}) follows.
\vskip0,1truecm
Analogously to (\ref{6}), we have
\begin{eqnarray}
\Vert\partial_x y(0,\cdot)\Vert_{L^2 (0,T)}^2 &\leq& C \left( \Vert y_0 \Vert^2 +\Vert\eta_0 \Vert_{L_\beta}^2 \right) +  \int_{0}^{T} \Vert y\Vert\Vert \sigma\Vert dt \nonumber\\
&{\le}& C \Vert \Theta_0\Vert_{\mathcal{H}}^2 +\left(\int_{0}^{T} \Vert \sigma\Vert dt\right)\displaystyle\max_{t \in [0,T]} \Vert y(\cdot,t)\Vert . \label{ux3}
\end{eqnarray}
Applying Young's inequality, \eqref{ux3} becomes
\begin{equation}
\Vert\partial_x y(0,\cdot)\Vert_{L^2 (0,T)}^2
\leq  C \left(\Vert \Theta_0\Vert_{\mathcal{H}}^2 +\displaystyle\max_{t \in [0,T]} \Vert y(\cdot,t)\Vert^2 \right)+ \frac{1}{2}\left( \int_{0}^{T} \Vert \sigma \Vert dt \right)^2 ,
\label{ux4}
\end{equation}
which together with (\ref{80}) yields
\begin{equation}
\Vert\partial_x y(0,\cdot)\Vert_{L^2 (0,T)}^2 \leq C \left(\Vert \Theta_0\Vert_{\mathcal{H}}^2 + \Vert \sigma\Vert_{L^1 ((0,T);L^2 (I))}^2 \right). \label{y022}
\end{equation}
A very similar argument as for the second estimate in \eqref{6} leads to
\begin{equation}\label{uu2}
\int_{0}^{T} \Vert \partial_x y \Vert^2 dt
\end{equation}
\begin{equation*}
\leq  C \left(  \displaystyle \max_{t \in [0,T]} \Vert y(\cdot,t)\Vert^2 +   \Vert\partial_x y(0,\cdot)\Vert_{L^2 ((0,T))}^2+ \Vert y_0 \|^2+ \Vert\eta_0 \Vert_{L_\beta}^2 + \Vert \sigma\Vert_{L^1 ((0,T);L^2 (I))}^2 \right),
\end{equation*}
which implies thanks to (\ref{80}) and (\ref{y022})
\begin{equation}
\int_{0}^{T} \Vert \partial_x y \Vert^2 dt \leq C \left(\Vert \Theta_0\Vert_{\mathcal{H}}^2 + \Vert \sigma\Vert_{L^1 ((0,T);L^2 (I))}^2 \right). \label{y02}
\end{equation}
Combining (\ref{80}) and (\ref{y02}), we have (\ref{81}).
\vskip0,1truecm
(ii) The proof of the second part of Theorem \ref{th2} is very similar to that of Proposition 4.1 in \cite{ro}.
\end{proof}

\subsection{Well-posedness of the problem (\ref{1})}\label{sub3} When
$\omega_3 =0$, the system \eqref{1} coincides with \eqref{lin}, and then the well-posedness of \eqref{1} is given in Theorem \ref{th1}. Therefore, in this subsection, we assume that $\omega_3 \ne 0$. Before stating the well-posedness result of \eqref{1}, let us first recall that the constant $C_1$ is defined in \eqref{81}. Moreover, let  $\gamma >0$ be the Sobolev embedding constant
\begin{equation}\label{injection}
\Vert v \Vert_{L^{\infty}(I)}^2 \leq \gamma \Vert v\Vert_{H^1(I)}^2 ,\quad v\in H^1 (I).
\end{equation}
\vskip0,1truecm
\begin{theorem}\label{th3}
Assume that ${\bf(H)}$ holds. Then, for every $\Theta_{0} =(y_0 ,\eta_0 )^T\in {\mathcal{H}}$ with
\begin{equation}\label{U0a}
\Vert \Theta_{0}\Vert_{\mathcal{H}}^2 < \frac{1}{4C_1\omega_3^2\gamma},
\end{equation}
the problem \eqref{UA} has a unique global solution $\Theta =(y ,\eta )^T$ satisfying the regularity \eqref{regularity0}, and consequently (\ref{1}) admits a unique global solution
$y\in\mathcal{M}$, which satisfies {the estimate}
\begin{equation}\label{yMCPhi0H}
\Vert y\Vert_{\mathcal{M}}\leq C\Vert \Theta_{0}\Vert_{\mathcal{H}},
\end{equation}
for some constant $C>0$.
\end{theorem}
\vskip0,1truecm
\begin{proof}
Let $\Theta_{0} \in {\mathcal{H}}$ satisfying (\ref{U0a}). Let $a>0$ such that
\begin{equation}\label{U0aplus}
\Vert \Theta_{0}\Vert_{\mathcal{H}}^2 \leq a <\frac{1}{4C_1\omega_3^2\gamma}.
\end{equation}
Next, consider the mapping
\begin{equation*}
\Gamma: \mathcal{M} \rightarrow \mathcal{M}
\end{equation*}
defined by $\Gamma({\tilde y})=y$, where $y$ is the solution of (\ref{fxt}) with the source term
\begin{equation*}
\sigma (x,t)=-\omega_3 {\tilde y}(x,t) \partial_x  {\tilde y}(x,t)
\end{equation*}
and initial data $\Theta_0$. In light of Theorem \ref{th2}, one can easily see that $\Gamma$ is well-defined and the following estimate holds:
\begin{equation*}
\Vert\Gamma({\tilde y})\Vert_{\mathcal{M}}^2 \le C_1 \left(\Vert \Theta_0\Vert_{\mathcal{H}}^2 +\omega_3^2\Vert {\tilde y} \partial_x {\tilde y} \Vert_{L^1 ((0,T);L^2 (I))}^2 \right).
\end{equation*}
Moreover, the embedding inequality (\ref{injection}) and the smallness condition (\ref{U0aplus}) on $\Vert \Theta_0\Vert_{\mathcal{H}}$ imply that
\begin{equation}
\Vert \Gamma({\tilde y}) \Vert_{\mathcal{M}}^2 \le C_1 \left( a + \omega_3^2\gamma \Vert {\tilde y} \Vert_{\mathcal{M}}^4 \right). \label{91}
\end{equation}
On the other hand, let ${\tilde y}_1 ,\,{\tilde y}_2\in \mathcal{M}$,
$y_1 =\Gamma({\tilde y}_1)$ corresponding to the initial data $\Theta_0$ and source term $\sigma_1 =-\omega_3 {\tilde y}_1 \partial_x  {\tilde y}_1$, $y_2 =\Gamma({\tilde y}_2)$ corresponding to the same initial data $\Theta_0$ and source term $\sigma_2 =-\omega_3 {\tilde y}_2 \partial_x {\tilde y}_2 $, and ${\hat y}=y_1 -y_2$. According to the definition of $\Gamma$, it is clear that $y$ is the solution of (\ref{fxt}) with the source term
\begin{equation*}
\sigma =-\omega_3 \left({\tilde y}_1 \partial_x  {\tilde y}_1 -{\tilde y}_2 \partial_x  {\tilde y}_2 \right)
\end{equation*}
and $(0,0)$ as initial data. Then, using (\ref{81}), we get
\begin{eqnarray*}
\Vert \Gamma({\tilde y}_1)-\Gamma({\tilde y}_2) \Vert_{\mathcal{M}}^2 &=& \Vert {\hat y} \Vert_{\mathcal{M}}^2\leq C_1 \omega_3^2 \Vert {\tilde y}_1 \partial_x  {\tilde y}_1 -{\tilde y}_2 \partial_x  {\tilde y}_2\Vert_{L^1 ((0,T);L^2 (I))}^2 \\
&\le& C_1 \omega_3^2\Vert \left({\tilde y}_1 -{\tilde y}_2 \right)\partial_x  {\tilde y}_2 +{\tilde y}_1\partial_x  \left({\tilde y}_1 -{\tilde y}_2 \right)\Vert_{L^1 ((0,T);L^2 (I))}^2 ;
\end{eqnarray*}
thus, exploiting Young's inequality and (\ref{injection}), we arrive at
\begin{equation}
\Vert \Gamma({\tilde y}_1)-\Gamma({\tilde y}_2) \Vert_{\mathcal{M}}^2 \le 2C_1 \omega_3^2\gamma  \left(  \Vert {\tilde y}_1 \Vert_{\mathcal{M}}^2 +\Vert {\tilde y}_2 \Vert_{\mathcal{M}}^2 \right) \Vert {\tilde y}_1 - {\tilde y}_2 \Vert_{\mathcal{M}}^2 .\label{92}
\end{equation}
\vskip0,1truecm
Now, we consider the restriction of $\Gamma$ to the closed ball
\begin{equation*}
B=\left\{ {\tilde y} \in \mathcal{M}, \, \Vert {\tilde y} \Vert_{\mathcal{M}}^2 \leq a \right \}.
\end{equation*}
This, together with (\ref{91}) and (\ref{92}), yields
\begin{equation}
\Vert\Gamma({\tilde y}) \Vert_{\mathcal{M}}^2 \leq C_1 \left( a + \omega_3^2\gamma a^2\right)\quad\hbox{and}\quad\Vert\Gamma({\tilde y}_1)-\Gamma({\tilde y}_2) \Vert_{\mathcal{M}}^2 \leq 4C_1\omega_3^2\gamma a\Vert {\tilde y}_1 -{\tilde y}_2 \Vert_{\mathcal{M}}^2 ,\quad {\tilde y}, \,{\tilde y}_1 ,\,{\tilde y}_2\in B.\label{93}
\end{equation}
Thereby, the map $\Gamma$ is well-defined and contractive on the ball
$B$ by virtue of (\ref{U0aplus}). Then, Banach Fixed Point Theorem leads to conclude that $\Gamma$ has a unique fixed element $y$, which turns out to be the unique solution to our problem (\ref{1}). Finally, since the energy $E$ of  (\ref{1}) is non-increasing, then the solution must be global and the estimate (\ref{yMCPhi0H}) can be obtained analogously to \eqref{81}.
\end{proof}

\section{Asymptotic behavior of the solutions to the KdVB equation}
\label{sec4}

Before announcing and proving our stability results, we consider, for a given $\Theta_0 \in\mathcal{H}$, the following additional hypothesis:
\begin{equation}\label{Ksc0a2a3}
\omega_0 >0\quad\hbox{or}\quad \left[\omega_0=0\,\,\hbox{and}\,\, c_0\vert\omega_2\vert+\frac{2}{3}\vert\omega_3\vert\gamma(1+c_0)\Vert\Theta_0\Vert_{\mathcal{H}} < 3\omega_1 \right],
\end{equation}
where $\gamma$ is defined in \eqref{injection} and $c_0$ is the smallest positive constant satisfying (Poincar\'e's inequality)
\begin{equation}\label{poincare}
\Vert v\Vert^2 \leq c_0 \Vert \partial_x v \Vert^2 ,\quad v\in H_0^1 (I).
\end{equation}
\vskip0,1truecm
\begin{theorem}\label{theorem1}
Assume that ${\bf(H)}$ and \eqref{Ksc0a2a3} hold, where $\Theta_0 \in\mathcal{H}$ satisfying \eqref{U0a} if $\omega_3\ne 0$. Then there exist positive constants $c$ and ${\tilde c}$ such that the solution $\Theta$ of \eqref{UA} satisfies the next two stability estimates:
\vskip0,1truecm
\begin{itemize}
\item[(i)] Case $\xi^{\prime}=0$:
\begin{equation}\label{decay1}
E (t)\leq {\tilde c} e^{-ct} , \quad t\in \mathbb{R}_+ .
\end{equation}
\vskip0,1truecm
\item[(ii)] Case $\xi^{\prime}\ne 0$:
\begin{equation}\label{decay2}
E (t)\leq {\tilde c} e^{-c\int_0^t \xi (s)ds}\left(1 +\int_0^t  \xi(s)e^{c\int_0^s \xi (\tau)d\tau}\int_s^{\infty} \beta (\tau)h(s,\tau)d\tau ds \right),\quad t\in \mathbb{R}_+ ,
\end{equation}
where, for $0\leq t\leq s$,
\begin{equation}\label{h}
h(t,s)=t+\left\vert\displaystyle \int_0^{s-t} y_1 (\tau)d\tau\right\vert^2 .
\end{equation}
\end{itemize}
\end{theorem}
\vskip0,1truecm
\begin{remark}\label{remark2}
When $\xi^{\prime}=0$ like \eqref{Example1} such that \eqref{Example11} holds, we get the exponential stability estimate \eqref{decay1} for \eqref{UA}. Nonetheless, when $\xi^{\prime}\ne 0$ like \eqref{Example2} such that \eqref{Example22} is satisfied, the decay rate of $E$ at infinity given by \eqref{decay2} depends on the ones of both $\beta$ and $y_1$. For example, let us consider the particular case $y_1\in L^{\infty}(\mathbb{R})$. Then, for some positive constant $C$,
\begin{equation*}
h(s,\tau)\leq C (\tau^2 +\tau),\quad 0\leq s\leq \tau.
\end{equation*}
Whereupon, integrating by parts, we find, for $t\in \mathbb{R}_+$,
\begin{eqnarray*}
\int_0^t  \xi(s) e^{c\int_0^s \xi (\tau)d\tau}\int_s^{\infty} h(s,\tau)\beta (\tau)d\tau ds &\leq& C \int_0^t  \partial_s \left(e^{c\int_0^s \xi (\tau)d\tau}\right)\int_s^{\infty} (\tau ^2 +\tau) \beta(\tau)d\tau ds\\
&\leq& C \left[ e^{c\int_0^t \xi (s)ds}\int_t^{\infty} (s^2 +s) \beta(s)ds \right.\\
&& \quad \quad \quad \left.+ \int_0^t e^{c\int_0^s \xi (\tau)d\tau} (s^2 +s)\beta (s)ds\right].
\end{eqnarray*}
Thus, \eqref{decay2} yields
\begin{eqnarray}\label{polydecaymodi}
E(t) &\leq& {\tilde c}e^{-c\int_0^t \xi (s)ds}+{\tilde c}C \int_t^{\infty} (s^2 +s) \beta (s)ds +{\tilde c}C e^{-c\int_0^t \xi (s)ds}\int_0^t e^{c\int_0^s \xi (\tau)d\tau} (s^2 +s)\beta (s)ds \nonumber\\
&\leq& {\tilde c}e^{-c\int_0^t \xi (s)ds}+2{\tilde c}C \int_t^{\infty} (s^2 +s) \beta (s)ds,
\end{eqnarray}
since $t\mapsto e^{-c\int_0^t \xi (s)ds}$ is non-increasing. If $\alpha$ is of the form \eqref{Example2} such that \eqref{Example22} holds and $d_1 >2$, then \eqref{polydecaymodi} leads to, for $b_2 ={\tilde c}\left(1+\frac{2d_1 d_2 C}{d_1 -2}\right)$ and $b_1 =\min\{cd_1 ,d_1 -2\}$,
\begin{equation*}
E(t)\leq b_2 (t+1)^{-b_1},\quad t\in \mathbb{R}_+ ,
\end{equation*}
since $\beta (s)=-\alpha^{\prime} (s)=d_1 d_2 (1+s)^{-(d_1 +1)}$ and $\xi (s)=d_1 (1+s)^{-1}$.
\end{remark}
\vskip0,1truecm
\begin{proof}
In order to prove \eqref{decay1} and \eqref{decay2}, and according to the assumption \eqref{Ksc0a2a3}, we distinguish the cases $\omega_0 >0$ and $\omega_0 =0$.

\subsection{Case 1: $\omega_0 >0$} Using \eqref{E}, \eqref{Eprime0}, \eqref{poincare} and the fact that $\beta^{\prime}\leq 0$, we have
\begin{equation}\label{Energy1}
E^{\prime} (t)\leq -\omega_0\Vert y_x\Vert^2\leq -\frac{\omega_0}{c_0}\Vert y\Vert^2 = -\frac{2\omega_0}{c_0} E (t)+\frac{\omega_0}{c_0} \int_0^{\infty} \beta (s) \eta^2 \,ds.
\end{equation}
Multiplying \eqref{Energy1} by $\xi$ and noticing that $E\geq 0$ and $\xi^{\prime}\leq 0$, we find
\begin{equation}\label{Energy3}
(\xi (t)E(t))^{\prime} \leq -\frac{2\omega_0}{c_0} \xi (t)E (t)+\frac{\omega_0}{c_0} \xi (t)\int_0^{\infty} \beta (s)\eta^2 \,ds.
\end{equation}
Now, we distinguish the two subcases $(i)$ and $(ii)$ considered in Theorem \ref{theorem1}.
\vskip0,1truecm
{\bf Subcase $1.1$: $\xi^{\prime} =0$}. Because $\xi$ is a positive constant, then, using \eqref{g01} and \eqref{Eprime0}, we deduce that
\begin{equation}\label{Energy30}
\xi (t) \int_0^{\infty} \beta (s) \eta^2 \,ds\leq -\int_0^{\infty} \beta^{\prime} (s)\eta^2 \,ds \leq -\frac{1}{\kappa} E^{\prime} (t).
\end{equation}
Therefore, combining \eqref{Energy3} and \eqref{Energy30}, we obtain, for the positive constant $c=\frac{2\omega_0\kappa\xi}{\omega_0+c_0\kappa\xi}$,
\begin{equation}\label{Energy31}
E^{\prime} (t)\leq -c E (t).
\end{equation}
Consequently, by integrating \eqref{Energy31}, we obtain \eqref{decay1} with ${\tilde c}=E(0)$.
\vskip0,1truecm
{\bf Subcase $1.2$: $\xi^{\prime} \ne 0$}. According to \eqref{Eprime0} and since $\beta^{\prime} \leq 0$, we have
\begin{equation*}
(\partial_x y)^2 (0,t)\leq -\frac{1}{\kappa} E^{\prime} (t),
\end{equation*}
which leads to
\begin{equation}\label{Energy4}
\int_0^{t} (\partial_x y)^2 (0,t) dt\leq \frac{1}{\kappa} E (0).
\end{equation}
On the other hand, applying Young's and H\"older's inequalities, we get, for
$0\leq t\leq s$,
\begin{eqnarray}\label{Energy50}
\left\vert\int_{t-s}^t  \partial_x y (0,\tau)d\tau\right\vert^2 &\leq& 2\left\vert\int_{t-s}^0 \partial_x y (0,\tau)d\tau\right\vert^2 +2\left\vert\int_{0}^t \partial_x y (0,\tau)d\tau\right\vert^2\\
&\leq& 2\left\vert\int_0^{s-t} y_1 (\tau)d\tau\right\vert^2 +2t\int_{0}^t (\partial_x y)^2 (0,\tau) d\tau\nonumber.
\end{eqnarray}
Thus, combining \eqref{Energy4} and \eqref{Energy50}, we get
\begin{equation}\label{Energy6}
\left\vert\int_{t-s}^t \partial_x y (0,\tau)d\tau\right\vert^2\leq c_1 h(t,s),\quad 0\leq t\leq s,
\end{equation}
for the positive constant $c_1=\max\left\{2,\frac{2}{\kappa}E(0)\right\}$ and $h$ is defined in \eqref{h}. Moreover, applying some arguments of \cite{chengues2, guesmia}, noticing that
$\xi (t)\leq \xi (s)$, for $0\leq s\leq t$, and using \eqref{g01} and \eqref{etaname}, we observe that
\begin{eqnarray}\label{BoundBound}
\xi (t)\int_0^{\infty} \beta (s)\eta^2 \,ds &=& \xi (t)\int_0^{t} \beta (s)\eta^2 \,ds +\xi (t)\int_t^{\infty} \beta (s)\eta^2 \,ds\\
&\leq& -\int_0^{t} \beta^{\prime} (s)\eta^2 \,ds + \xi (t)\int_t^{\infty} \beta (s)\left\vert \int_{t-s}^t \partial_x y (0,\tau)d\tau\right\vert^2 \,ds \nonumber\\
&\leq& -\int_0^{\infty} \beta^{\prime} (s)\eta^2 \,ds + \xi (t)\int_t^{\infty} \beta (s)\left\vert \int_{t-s}^t \partial_x y (0,\tau)d\tau\right\vert^2 \,ds. \nonumber
\end{eqnarray}
Consequently, using \eqref{Eprime0} and \eqref{Energy6}, we deduce from \eqref{BoundBound} that
\begin{equation}\label{Energy7}
\xi (t)\int_0^{\infty} g (s)\eta^2 \,ds\leq - \frac{1}{\kappa} E^{\prime} (t) +c_1 \xi (t)\int_t^{\infty} \beta (s)h(t,s) \,ds.
\end{equation}
We set
\begin{equation}\label{Fxi}
F(t)=\left(\xi (t)+\frac{\omega_0}{c_0\kappa}\right)E(t).
\end{equation}
Because $\xi\geq 0$ and $\xi^{\prime}\leq 0$, we see that
\begin{equation}\label{equivalence}
\frac{\omega_0}{c_0\kappa} E (t) \leq F(t) \leq \left(\xi (0) +\frac{\omega_0}{c_0\kappa}\right) E(t).
\end{equation}
Exploiting \eqref{Energy3}, \eqref{Energy7} and the right inequality in \eqref{equivalence}, we obtain
\begin{equation}\label{Energy80}
F^{\prime} (t)\leq -c\xi (t)F (t)+\frac{c_1 \omega_0}{c_0}\xi (t)\int_t^{\infty} \beta (s)h(t,s) \,ds,
\end{equation}
for the positive constant $c=\frac{2\omega_0\kappa}{\omega_0+c_0\kappa\xi (0)}$, this implies that
\begin{equation}\label{Energy8}
\partial_t\left[e^{c\int_0^t \xi (s)ds}F(t)-\frac{c_1 \omega_0}{c_0}\int_0^{t}\xi (s)e^{c\int_0^s \xi (\tau)d\tau}\int_s^{\infty} \beta (\tau)h(s,\tau) \,d\tau ds\right] \leq 0.
\end{equation}
Integrating \eqref{Energy8}, we find
\begin{equation}\label{Energy9}
F(t)\leq e^{-c\int_0^t \xi (s)ds} \left[F(0) +\frac{c_1 \omega_0}{c_0} \int_0^{t} \xi (s)e^{c\int_0^s \xi (\tau)d\tau}\int_s^{\infty} \beta (\tau)h(s,\tau) \,d\tau ds\right],
\end{equation}
and therefore, using \eqref{equivalence}, we reach \eqref{decay2} with ${\tilde c}=\frac{c_0\kappa}{\omega_0}\left\{\left(\xi (0)+\frac{\omega_0}{c_0\kappa}\right)E(0), \frac{c_1\omega_0}{c_0}\right\}$.

\subsection{Case 2: $\omega_0 =0$} Similarly to \eqref{ux} and \eqref{px}, multiplying (\ref{1})$_1$ by $xy(x,t)$, integrating by parts over $I$ and using the boundary conditions in \eqref{1} and \eqref{etatheta}, we obtain
\begin{eqnarray*}
3\omega_1 \int_0^1 \left( \partial_x y \right)^2 dx &=& \omega_2 \int_0^1 y^2 dx +\frac{2\omega_3}{3} \int_0^1 y^3 dx -\partial_t\left(\displaystyle\int_0^1 xy^2\,dx \right)\\
&& +\omega_1\left( \int_0^{\infty} \beta(s)  \eta ds \right)^2 + 2 \omega_1\omega_4 \partial_x y(0,t)  \int_0^{\infty} \beta(s)  \eta ds +\omega_1 \omega_4^2 \left( \partial_x y \right)^2 (0,t) dt .
\end{eqnarray*}
Using Young's inequality, \eqref{di2} and \eqref{Eprime0}, we see that, for some positive constant $c_1$,
\begin{equation*}
\omega_1\left( \int_0^{\infty} \beta(s)  \eta ds \right)^2 + 2 \omega_1\omega_4 \partial_x y(0,t)  \int_0^{\infty} \beta(s)  \eta ds +\omega_1 \omega_4^2 \left( \partial_x y \right)^2 (0,t) dt
\end{equation*}
\begin{eqnarray*}
&\leq& c_1\left[\left( \partial_x y \right)^2 (0,t)-\int_0^{\infty} \beta^{\prime}(s)\eta^2 ds \right]\\
&\leq& -\frac{c_1}{\kappa}E^{\prime} (t).
\end{eqnarray*}
Therefore, by combining the above two formulas and using \eqref{poincare}, we arrive at
\begin{equation}\label{xy1}
\left(3\omega_1-c_0\vert\omega_2\vert\right)\Vert \partial_x y\Vert^2 \leq \frac{2\omega_3}{3}\int_0^1 y^3 dx-\partial_t \left(\int_0^1 xy^2 dx+\frac{c_1}{\kappa}E(t)\right) .
\end{equation}
On the other hand, using \eqref{E}, \eqref{injection}, \eqref{poincare} and
H\"older's inequality and noticing that $E^{\prime}\leq 0$, we see that
\begin{eqnarray}\label{xy2}
\left\vert\int_0^1 y^3 dx\right\vert &\leq& \Vert y\Vert_{L^{\infty}(I)}^2 \int_0^1 \vert y\vert dx\leq \gamma \left(\Vert y\Vert^2 +\Vert\partial_x y\Vert^2\right) \Vert y\Vert \\
&\leq& \gamma (1+c_0)\Vert\partial_x y\Vert^2 {\sqrt{2E(t)}} \leq \gamma (1+c_0){\sqrt{2E(0)}}\Vert\partial_x y\Vert^2.\nonumber
\end{eqnarray}
Thus, by combining \eqref{xy1} and \eqref{xy2}, it follows that
\begin{equation*}
\left[3\omega_1-c_0\vert\omega_2\vert-\frac{2}{3}\vert\omega_3\vert\gamma(1+c_0)\Vert\Theta_0\Vert_{\mathcal{H}}\right]\Vert \partial_x y\Vert^2 \leq -\partial_t \left(\int_0^1 xy^2 dx+\frac{c_1}{\kappa}E(t)\right).
\end{equation*}
Consequently, combining the latter with \eqref{Ksc0a2a3}, we have
\begin{equation}\label{xy4}
\Vert \partial_x y\Vert^2 \leq -C_0\partial_t \left(\int_0^1 xy^2 dx+\frac{c_1}{\kappa}E(t)\right) ,
\end{equation}
where
$$C_0 =\frac{1}{3\omega_1-c_0\vert\omega_2\vert
-\frac{2}{3}\vert\omega_3\vert\gamma(1+c_0)\Vert\Theta_0\Vert_{\mathcal{H}}}.
$$
Hence, we deduce from \eqref{E}, \eqref{poincare} and \eqref{xy4} that, for $C_1 =\frac{c_0 C_0}{2}$,
\begin{equation}\label{MainInequ}
E(t)\leq \frac{c_0}{2}\Vert\partial_x y\Vert^2 +\frac{1}{2}\int_0^{\infty} \beta(s) \eta^2 ds\leq -C_1\partial_t \left(\int_0^1 xy^2 dx+\frac{c_1}{\kappa}E(t)\right)+\frac{1}{2}\int_0^{\infty} \beta(s) \eta^2 ds.
\end{equation}
\vskip0,1truecm
{\bf Subcase $2.1$: $\xi^{\prime}=0$}. Because $\xi$ is a positive constant, then multiplying \eqref{MainInequ} by $\xi$ and exploiting \eqref{Eprime0} and the right inequality in \eqref{g01}, we get
\begin{equation}\label{MainInequB}
\partial_t\left[\frac{1+2c_1 C_1\xi}{2\kappa}E(t)+C_1\xi\int_0^1 xy^2 dx\right]\leq -\xi E(t).
\end{equation}
Let us consider the function
\begin{equation*}
F(t)=\frac{1+2c_1 C_1\xi}{2\kappa}E(t)+C_1\xi\int_0^1 xy^2 dx.
\end{equation*}
We see that
\begin{equation}\label{MainInequC}
\frac{1+2c_1 C_1\xi}{2\kappa}E(t)\leq F(t)\leq \left(\frac{1+2c_1 C_1\xi}{2\kappa}+2C_1\xi\right)E(t),
\end{equation}
thus, using \eqref{MainInequB} and the right inequality in \eqref{MainInequC}, we find, for $c=\frac{2\kappa\xi}{1+2c_1 C_1\xi+4C_1\kappa\xi}$, that $F^{\prime}\leq -cF$, which, by integrating, implies that
\begin{equation*}
F(t)\leq F(0)e^{-ct}.
\end{equation*}
Hence, according to \eqref{MainInequC}, we deduce that \eqref{decay1} is satisfied with ${\tilde c}=\frac{1+2c_1 C_1\xi+4C_1\kappa\xi}{1+2c_1 C_1\xi}$.
\vskip0,1truecm
{\bf Subcase $2.2$: $\xi^{\prime}\ne 0$}. Multiplying \eqref{MainInequ} by
$\xi (t)$ and exploiting \eqref{Energy7}, we find
\begin{equation}\label{MainInequA}
\frac{1}{2\kappa}E^{\prime}(t)+C_1\xi (t)\partial_t \left(\int_0^1 xy^2 dx+\frac{c_1}{\kappa}E(t)\right)\leq -\xi (t)E(t)+\frac{c_1}{2}\xi (t)\int_t^{\infty} \beta(s) h(t,s)ds.
\end{equation}
Subsequently, let
\begin{equation}
F(t)= \frac{1}{2\kappa} E(t) +C_1 \xi (t)\left(\int_0^1 xy^2 dx+\frac{c_1}{\kappa}E(t)\right). \label{gt}
\end{equation}
Since $\xi^{\prime}\leq 0$ and
\begin{equation*}
0\leq \xi(t)\left(\int_0^1 xy^2 dx+\frac{c_1}{\kappa}E(t)\right)\leq \xi(0)\left( 2+\frac{c_1}{\kappa}\right)E(t),
\end{equation*}
it follows that
\begin{equation}\label{equivalence0}
\frac{1}{2\kappa} E(t) \leq F(t) \leq \left[\frac{1}{2\kappa}+\xi(0)\left( 2+\frac{c_1}{\kappa}\right)\right]E(t).
\end{equation}
Then, using \eqref{MainInequA}, \eqref{gt} and the right inequality in \eqref{equivalence0}, and noticing again that $\xi^{\prime}\leq 0$, we obtain, for $c=\frac{2\kappa}{1+2\xi (0)(2\kappa +c_1)}$,
\begin{equation}\label{MainInequ1}
F^{\prime}(t)\leq -c\xi(t)F(t)+\frac{c_1}{2}\xi (t)\int_t^{\infty} \beta(s) h(t,s) ds,
\end{equation}
which is similar to \eqref{Energy80}, and hence the proof of \eqref{decay2} can be achieved as in the previous subcase 1.2.
\end{proof}
\vskip0,1truecm
\begin{remark}
It is interesting to mention that the well-posedness and stability results shown for the KdVB equation include the case $\omega_0=0$. This means that our findings remain valid for the KdV equation.
\end{remark}

\section{Application to the Kuramoto-Sivashinsky equation}\label{sec5}

In this section, we extend our results to the well-known fourth-order KS equation with boundary infinite memory
\begin{equation}
\left\{
\begin{array}{ll}
{\partial_{t} y (x,t)} +\nu_0 \partial_x^{4} y{(x,t)} +\nu_1 y{(x,t)} \partial_x y{(x,t)} +\nu_2 \partial_x^2 y{(x,t)}=0, &  x \in  I, \, t >0,\\
y (0,t) =y(1,t) =\partial_x^2 y(1,t) =0, & t \in\mathbb{R}_+ ,\\
\partial_x^2 y(0,t)=\nu_3 \partial_x y(0,t)+\displaystyle\int_0^{\infty} \alpha (s) \partial_x y(0,t-s)  ds, &t \in\mathbb{R}_+ , \\
y(x,0) =y_{0} (x) , & x \in I,\\
\partial_x y(0,-t)=y_1 (t), & t\in\mathbb{R}_+ 
\end{array}
\right.
\label{11}
\end{equation}
in which $\nu_i$ are real constants (physical parameters) satisfying the following hypothesis ${\bf ({\tilde H})}$:
\vskip0,1truecm
\begin{itemize}
\item The constants $\nu_0, \nu_2$ and $\nu_3$ satisfy
\begin{equation}
\nu_0 >0, \quad \nu_3 > 0 \quad\hbox {and}\quad 0<\nu_2 < \pi^2 \nu_0.\label{ai1}
\end{equation}
\item The memory kernel $\alpha$ satisfies \eqref{f}, \eqref{fi} and \eqref{xi}. In turn, instead of \eqref{f2}, $\alpha$ obeys the condition
\begin{equation}\label{f3}
\int_0^{\infty}\frac{-\alpha^{\prime}(s)}{\xi (s)}dx:=\alpha_0 <
\frac{2\nu_0\nu_3}{\left\vert 1-\nu_0\right\vert^2}\quad \hbox {if}\,\,\nu_0\ne 1.
\end{equation}
\end{itemize}
No condition is considered on $\alpha_0$ if $\nu_0 = 1$.  
\vskip0,1truecm
The reader who is interested in a literature review of the KS equation can consult  \cite{chen1, chen2} and the references therein.
\vskip0,1truecm
Next, we shall adopt the same notations (\ref{etaname}) and (\ref{Lg})-(\ref{HInner}) as in Section \ref{sec2}, so (\ref{etadef}) and (\ref{etatheta}) are valid, and accordingly the problem (\ref{11}) can be formulated in $\mathcal{H}$ as follows:
\begin{align}\label{UB}
\begin{cases}
\partial_t \Lambda (t) = \mathcal{Q} \Lambda(t),\quad t>0,\\
\Lambda(0) = \Lambda_{0},
\end{cases}
\end{align}
where $\Lambda=(y,\eta)^T$, $\Lambda_{0}=(y_{0},\eta_0)^T$ and $\mathcal{Q}$  is the nonlinear operator defined by
\begin{equation*}
\left\{
\begin{array}{l}
\mathcal{D} (\mathcal{Q})=\left\{
\begin{array}{l}
\Lambda \in \mathcal{H}; \, \mathcal{Q} \Lambda \in \mathcal{H} ,\, y (0,\cdot)=y(1,\cdot)=\partial_x^2 y(1,\cdot)=0 ,\,\eta (\cdot,0)=0, \\[1mm]
\partial_x^2 y(0,\cdot)= \nu_3 \partial_x y(0,\cdot)+\displaystyle\int_0^{\infty} \beta (s) \eta (\cdot,s)  ds
\end{array}\right\} \\[3mm]
\mathcal{Q} \Lambda =\left(
\begin{array}{c}
-\nu_0  \partial_x^4 y -\nu_1 y  \partial_x y -\nu_2  \partial_x^2 y\\
\partial_x y(0,\cdot)-\partial_s \eta
\end{array}
\right).
\end{array}
\right.
\end{equation*}
Let us note that the spaces $H_0^1 (I)$ and $H_0^1 (I) \cap H^2 (I)$ will be, respectively, equipped with the equivalent norms $ \| \partial_x \cdot \|$ and $ \| \partial_x^2 \cdot \|$ in light of the following Wirtinger's inequalities \cite{H, wang}:
\begin{eqnarray}
\displaystyle \pi^2 \int_{0}^{1} u^2 (x) \, dx  &\leq& \int_{0}^{1} (\partial_x u)^2 (x) \, dx, \; \forall u \in H_0^1 (I); \label{wi1} \\
\displaystyle \pi^2 \int_{0}^{1}  (\partial_x u)^2 (x) \, dx  &\leq&    \int_{0}^{1} (\partial_x^2 u)^2 (x) \, dx, \; \forall u \in H_0^1 (I) \cap H^2 (I). \label{wi2}
\end{eqnarray}
Moreover, for $T>0$, we consider the space
$$
\mathcal{S}=C \left( [0,T]; \, L^2 (I) \right) \cap L^2  \left( (0,T); \, H^1_0 (I)\cap H^2 (I) \right), 
$$
whose norm is
$$
\| \cdot \|_{\mathcal{S}}^2=\| \cdot \|_{C ( [0,T]; \, L^2 (I) )}^2 + \| \cdot \|_{L^2  ( (0,T); \, H^2 (I))}^2.
$$
Thereafter, we merely argue as for the KdVB equation. 

\subsection{The linearized system associated to (\ref{11})} The linear system associated to (\ref{11}) (that is (\ref{11}) with $\nu_1=0$) can be written in $\mathcal{H}$ as follows:
\begin{equation}
\left\{
\begin{array}{ll}
\partial_t {\Lambda} (t)=\mathcal{K} \Lambda(t),\quad t>0,\\
\Lambda(0)=\Lambda_{0},
\end{array}
\right. \label{siks}
\end{equation}
where $\mathcal{K}$ is the linear operator defined by
\begin{equation}
\left\{
\begin{array}{l}
{\mathcal{D}}\left(\mathcal{K}\right)=\left\{
\begin{array}{l}
\Lambda \in \mathcal{H};\, \mathcal{K} \Lambda \in \mathcal{H},\,
 y (0,\cdot)=y(1,\cdot)=\partial_x^2 y(1,\cdot)=0 ,\,\eta (\cdot,0)=0, \\[1mm] \partial_x^2 y(0,\cdot)= \nu_3 \partial_x y(0,\cdot)+\displaystyle\int_0^{\infty} \beta (s) \eta (\cdot,s)  ds
 \end{array} \right\},\\[3mm]
\mathcal{K} \Lambda=
\begin{array}{l}
 \left(
\begin{array}{c}
-\nu_0 \partial_x^4 y-\nu_2 \partial_{x}^2 y \\
\partial_x y(0,\cdot)-\partial_s \eta
\end{array}
\right)
\end{array}.
\end{array}
\right.
\label{1.62ks}
\end{equation}
\vskip0,1truecm
\begin{theorem}\label{th510}
Assume that ${\bf ({\tilde H})}$ hold. Then we have:
\begin{itemize}
\item[(i)] The linear operator $\mathcal{K}$ defined by (\ref{1.62ks}) generates a $C_{0}$-semigroup of contractions $e^{t\mathcal{K}}$ and hence, given an initial data $\Lambda_{0}\in{\mathcal{D}}(\mathcal{K})$, the Cauchy problem (\ref{siks}) has a unique classical solution
\begin{equation}
\Lambda  \in C(\mathbb{R}_+;{\mathcal{D}}(\mathcal{K})) \cap C^{1}(\mathbb{R}_+;\mathcal{H}). \label{regularityks}
\end{equation}
If $\Lambda_{0}\in \mathcal{H}$, then (\ref{siks}) has a mild solution
\begin{equation}
\Lambda \in  C(\mathbb{R}_+;\mathcal{H}). \label{regularity0ks}
\end{equation}
\vskip0,1truecm
\item[(ii)] For each $\Lambda_{0}\in \mathcal{H}$ and $T>0$, there exists a positive constant $C$ such that the solution $\Lambda$ of \eqref{siks} stemmed from the initial data $\Lambda_0$ satisfies the first estimate of \eqref{6} (with $\Lambda_{0}$ instead of $\Theta_{0}$) as well as
\begin{equation}
\Vert \partial_x^2 y \Vert_{L^2 ((0,T);L^2 (I) )}^2  \leq  C  \Vert \Lambda_0\Vert_{\mathcal{H}}^2 . \label{6ks}
\end{equation}
\vskip0,1truecm
Lastly, the mapping
\begin{equation}
\Upsilon: \Lambda_{0} \in \mathcal{H} \mapsto \Upsilon \left( \Lambda_{0}\right)=\Lambda (\cdot) = e^{\cdot \mathcal{K}} \Lambda_{0} \in \mathcal{S} \times C \left( [0,T]; \, L_{\beta} \right)
\end{equation}
is continuous.
\end{itemize}
\label{th1ks}
\end{theorem}
\vskip0,1truecm
\begin{proof}
Arguing as for \eqref{di1} and using (\ref{1.62ks}), we have that, for any 
$\Lambda=(y,\eta)^T$ in ${\mathcal{D}}({\mathcal{K}})$,
\begin{eqnarray*}
\langle{\mathcal{K}}\Lambda,\Lambda \rangle_{\mathcal{H}}
&=& -\nu_0 \nu_3 (\partial_{x}y)^2 (0,\cdot)-\nu_0 \Vert \partial_x^2 y \Vert^2 +\nu_2 \Vert \partial_x y \Vert^2 + (1-\nu_0)\partial_x y(0,\cdot)\int_0^{\infty} \beta(s)\eta ds\\
&& +\frac{1}{2}\int_0^{\infty}\beta^{\prime}(s)\eta^2 ds.
\end{eqnarray*}
In view of \eqref{di2}, \eqref{di3} and \eqref{wi2}, the latter gives
\begin{eqnarray}
\langle{\mathcal{K}}\Lambda,\Lambda \rangle_{\mathcal{H}}
&\le& \left[\vert 1-\nu_0 \vert\dfrac{\epsilon}{2}-\nu_0 \nu_3\right] (\partial_{x}y)^2 (0,\cdot)+\frac{1}{2}\left[ 1-\dfrac{\alpha_0}{\epsilon}\vert 1-\nu_0\vert \right] \int_0^{\infty}\beta^{\prime}(s)\eta^2 ds \nonumber \\
&& +\left(\dfrac{\nu_2}{\pi^2}-\nu_0 \right) \Vert \partial_x^2 y \Vert^2, \label{diks}
\end{eqnarray}
for any $\epsilon >0$. By virtue of ${\bf ({\tilde H})}$, one can choose 
$\epsilon$ as follows
\begin{equation}\label{eps}
\alpha_0 \vert 1-\nu_0\vert < \epsilon < \dfrac{2 \nu_0 \nu_3}{\vert 1-\nu_0\vert}\quad\hbox{if}\,\,\nu_0\ne 1 ,
\end{equation}
and consequently \eqref{diks} leads to (in both cases $\nu_0 = 1$ and 
$\nu_0\ne 1$)
\begin{equation}\label{K}
\langle{\mathcal{K}}\Lambda,\Lambda \rangle_{\mathcal{H}}
\le -\vartheta (\partial_{x}y)^2 (0,\cdot) -\vartheta\Vert \partial_x^2 y \Vert^2 +\vartheta \int_0^{\infty}\beta^{\prime}(s)\eta^2 ds,
\end{equation}
where
\begin{equation}\label{vart}
\vartheta=\min \left\{\frac{1}{2}\left[ 1-\dfrac{\alpha_0}{\epsilon}\vert 1-\nu_0\vert \right], \nu_0 \nu_3-\vert 1-\nu_0 \vert\dfrac{\epsilon}{2} ,\nu_0 -\frac{\nu_2}{\pi^2} \right\}.
\end{equation}
Clearly, $\vartheta$ is a well-defined positive number in view of ${\bf ({\tilde H})}$. This, together with \eqref{ai1} and \eqref{K}, implies that $\mathcal{K}$ is dissipative.
\vskip0,1truecm
Next, we show that $\lambda I-\mathcal{K}$ is onto $\mathcal{H}$, for any 
$\lambda>0$. Indeed, given $(z,f)^T$ in $\mathcal{H}$, we seek $(y,\eta)^T$ in ${\mathcal{D}}({\mathcal{K}})$ so that
\begin{equation}\label{onto}
\left\{
\begin{array}{l}
    \lambda y +\nu_0 \partial_x^4 y+\nu_2 \partial_x^2 y =z,\\
\lambda \eta -\partial_x y(0,\cdot)+\partial_s \eta =f, \\
    y (0,\cdot)=y(1,\cdot)=\partial_x^2 y(1,\cdot)=0,  \\
    \eta (\cdot,0)=0, \\
    \partial_x^2 y(0,\cdot)= \nu_3 \partial_x y(0,\cdot)+\displaystyle\int_0^{\infty} \beta (s) \eta (\cdot,s)  ds.
  \end{array}
  \right.
\end{equation}
Solving \eqref{onto}$_2$ and using \eqref{onto}$_4$, we obtain
\[ \eta(\cdot,s)=\displaystyle \int_{0}^{s} e^{-\lambda (s-r)}  \left( \partial_x y(0,\cdot) +f(r) \right) dr.\]
Thereby, it amounts to solving the following problem:
\begin{equation}\label{yp}
  \left\{
  \begin{array}{l}
    \lambda y +\nu_0 \partial_x^4 y+\nu_2 \partial_x^2 y =z,\\
    y (0,\cdot)=y(1,\cdot)=\partial_x^2 y(1,\cdot)=0,  \\
    \partial_x^2 y(0,\cdot)=\left( \nu_3 +\displaystyle \int_0^{\infty} \int_{0}^{s} \beta (s) e^{-\lambda (s-r)}  ~ dr  ds \right) \partial_x y(0,\cdot) +
     \displaystyle \int_0^{\infty} \int_{0}^{s} \beta (s) e^{-\lambda (s-r)} f(r) ~ dr  ds,
  \end{array}
  \right.
\end{equation}
which has the weak formulation
\begin{align*}
&  \displaystyle \int_{0}^{1} \left( \lambda y \phi + \nu_0 \partial_x^2 y \partial_x^2 \phi -\nu_2 \partial_x y  \partial_x \phi \right) dx +
 \nu_0 \left( \nu_3 +\displaystyle \int_0^{\infty} \int_{0}^{s} \beta (s) e^{-\lambda (s-r)}  ~ dr  ds \right) \partial_x y(0) \partial_x \phi (0)     \\
 &   =-\nu_0\partial_x \phi (0) \displaystyle \int_0^{\infty} \int_{0}^{s} \beta (s) e^{-\lambda (s-r)} f(r) ~ dr  ds +\int_{0}^{1} z \phi dx,
\end{align*}
for any $\phi$ in $H_0^1(I)\cap H^2 (I)$. Lastly, Lax-Milgram Theorem (see for instance \cite{br}) permits to claim the existence and uniqueness of a solution $y$ in $H_0^1(I)\cap H^2 (I)$ to the last problem and then by virtue of standard arguments used for elliptic linear equations, we can check that  $y \in H_0^1(I) \cap H^4(I) $ and recover the boundary conditions. Thus, the operator $\lambda I-{\mathcal{K}}$ is onto $\mathcal{H}$. The assertions (i) immediately follow from the fact that $\mathcal{K}$ is  a closed and densely defined operator and the semigroups theory of linear operators \cite{Pa}.
\vskip0,1truecm
With regard to the second item of the theorem, it suffices to establish it for solutions of \eqref{siks} stemmed from the domain $\mathcal{D}(\mathcal{K})$ by means of use a standard argument of density. Then, multiply the first (resp. second) equation of \eqref{siks} by $y$ (resp. $\beta(s) \eta$) and then integrate over $I$ (resp. $[0,\infty)$), we get (after performing similar computations as for \eqref{K})
\[
\vartheta (\partial_{x}y)^2 (0,\cdot) + \vartheta\Vert \partial_x^2 y \Vert^2 -\vartheta\displaystyle\int_0^{\infty} \beta^{\prime} (s) \eta (\cdot,s) ds \le -\left\langle\mathcal{K} \Lambda,\Lambda \right\rangle_{\mathcal{H}}=-\left\langle\partial_t\Lambda,\Lambda \right\rangle_{\mathcal{H}}=-\partial_t\left(\frac{1}{2}\Vert\Lambda\Vert_{\mathcal{H}}^2\right),
\]
where $\vartheta$ is defined by \eqref{vart}, therefore, by integrating over $[0,T]$, it follows that
\[
\vartheta \|(\partial_{x}y) (0,\cdot)\|_{L^2(0,T)}^2 + \vartheta\Vert \partial_x^2 y \Vert^2_{L^2\left( (0,T); L^2(I)\right)} -\vartheta\int_0^T \displaystyle\int_0^{\infty} \beta^{\prime} (s) \eta (t,s) dsdt\le \frac{1}{2}\| \Lambda_0 \|_{\mathcal{H}}^2.
\]
This together with the contraction of the semigroup $e^{t \mathcal{K}}$ leads to the desired results.
\end{proof}

\subsection{A non-homogeneous linear system associated to (\ref{11})} The next step is to consider the linear system (\ref{siks}) but with a source term $z: (x,t)\in I\times\mathbb{R}_+\mapsto z(x,t)\in \mathbb{R}$, namely,
\begin{equation}
\left\{
\begin{array}{ll}
\partial_t {\Lambda} (t)=\mathcal{K} \Lambda(t)+(z(x,t),0),\quad t>0,\\
\Lambda(0)=\Lambda_{0},
\end{array}
\right. \label{fxtks}
\end{equation}
whose energy is also defined by \eqref{E} (with $\Lambda$ instead of 
$\Theta$). We have the following result:
\vskip0,1truecm
\begin{theorem} Assume that ${\bf ({\tilde H})}$ holds. Then we have:
\begin{itemize}
\item[(i)] If $\Lambda_{0}=(y_{0},\eta_0 )^T \in \mathcal{H}$ and
$z \in L^1 ((0,T);L^2 (I))$, then there exists a unique mild solution
$\Lambda=(y,\eta^t)^T$ of (\ref{fxtks}) such that
$$
\Lambda \in \mathcal{S} \times C([0,T]; L_{\beta} ),
$$
and
\begin{equation}
\left\{
\begin{array}{l}
\Vert \Lambda \Vert_{C([0,T]; \mathcal{H})}^2  \leq  C_0 \left(\Vert \Lambda_{0}\Vert_{\mathcal{H}}^2 +\Vert z \Vert_{L^1 ((0,T);L^2 (I))}^2\right),\\
\Vert \partial_x y (0,\cdot) \Vert_{L^2 (0,T)}^2  \leq  C_1 \left(\Vert \Lambda_{0}\Vert_{\mathcal{H}}^2 +\Vert z \Vert_{L^1 ((0,T);L^2 (I))}^2\right),\\
\Vert y \Vert_{\mathcal{S}}^2  \leq  C_2 \left(  \Vert \Lambda_{0}\Vert_{\mathcal{H}} ^2 + \Vert z \Vert_{L^1 ((0,T);L^2 (I))}^2 \right),\label{81ks}
\end{array}
\right.
\end{equation}
for some positive constants $C_0 ,\,C_1>0$ and $C_2$ independent of 
$\Lambda_{0}$ and $z$.
\vskip0,1truecm
\item[(ii)] Given $y \in \mathcal{S}$, we have $y \partial_x y \in L^1 ((0,T);L^2 (I))$ and the map
\begin{equation*}
\Delta: y \in \mathcal{S} \mapsto y \partial_x y \in L^1 ((0,T);L^2 (I))
\end{equation*}
is continuous.
\end{itemize}
\label{th2ks}
\end{theorem}
\vskip0,1truecm
\begin{proof} The contraction of the semigroup $e^{t\mathcal{K}}$ and the fact that $z \in L^1 ((0,T);L^2 (I) )$ allow us to deduce the existence, uniqueness and smoothness results \cite{Pa}. For the estimates \eqref{81ks}, using once again a density argument and recalling that the energy of (\ref{fxtks}) is defined by \eqref{E}, we can obtain as for \eqref{K}
\begin{equation}
E^{\prime}(t) \leq -\vartheta (\partial_{x}y)^2 (0,\cdot) -\vartheta\Vert \partial_x^2 y(t) \Vert^2 +\vartheta\displaystyle\int_0^{\infty} \beta^{\prime} (s) \eta (\cdot,s) ds+ \int_{0}^{1} y(x,t)  z(x,t) dx.
\label{fi4}
\end{equation}
Subsequently, we integrate \eqref{fi4} and then use Cauchy-Schwarz and Young's inequalities to reach
\begin{eqnarray}
&& \Vert \Lambda\Vert_{\mathcal{H}}^2 + 2\vartheta \Vert \partial_x y (0,\cdot) \Vert_{L^2 (0,T)}^2 +2\vartheta\Vert \partial_x^2 y \Vert_{L^2 ((0,T);L^2 (I) )}^2 -2\vartheta\displaystyle\int_0^{\infty} \beta^{\prime} (s) \eta (\cdot,s) ds \nonumber\\
&& \quad \leq \Vert \Lambda_0\Vert_{\mathcal{H}}^2  + \delta \displaystyle\max_{t\in[0,T]}\Vert y(\cdot,t)\Vert^2 ds +\dfrac{1}{\delta}  \left(\int_{0}^{T} \Vert z (\cdot,t)\Vert dt\right)^2, \label{83ks}
\end{eqnarray}
for any $\delta >0$. Finally, invoking \eqref{wi1} and \eqref{wi2} and picking up $\delta$ small enough, we obtain the desired estimates \eqref{81ks}.
\vskip0,1truecm
Concerning the proof the second part (ii), the reader is referred to \cite{chengues2}.
\end{proof}

\subsection{Well-posedness of (\ref{11})} Based on the above discussion, one can obtain analogously to Theorem \ref{th3} (see also \cite[Theorem 2.4]{chengues2}) the following theorem:
\vskip0,1truecm
\begin{theorem}\label{kssol}
Assume that ${\bf ({\tilde H})}$ holds. Given $T>0$, there exist two positive constants $M_0$ and $M$ such that for every initial condition $\Lambda_0 =(y_0,\eta_0)^T \in \mathcal{H}$ satisfying 
\begin{equation}
\|\Lambda_0\|_{\mathcal{H}} \le M_0, \label{SCKSid}
\end{equation}
the problem \eqref{11} has a unique solution $y \in \mathcal{S}$. Moreover, we have 
\begin{equation*}
\|y\|_{\mathcal{S}} \le M \|\Lambda_0\|_{\mathcal{H}}.
\end{equation*}
\end{theorem}

\subsection{Stability of (\ref{11})} Our stability results for (\ref{UB}) are the same as for \eqref{UA}, more precisely, we have the next theorem.  
\vskip0,1truecm
\begin{theorem}\label{theorem5511}
Assume that ${\bf({\tilde H})}$ holds and $\Lambda_0 \in\mathcal{H}$ satisfying \eqref{SCKSid}. Then there exist positive constants $c$ and 
${\tilde c}$ such that the solution $\Lambda$ of \eqref{UB} satisfies the stability estimates \eqref{decay1} and \eqref{decay2}, where $h$ is defined in \eqref{h}.
\end{theorem}
\vskip0,1truecm
\begin{proof}
Because \eqref{fxtks} with $z=-\nu_1 y\partial_x y$ is reduced to (\ref{UB}), then, from \eqref{fi4} with $z=-\nu_1 y\partial_x y$, we conclude that 
\begin{equation*}
E^{\prime}(t) \leq -\vartheta (\partial_{x}y)^2 (0,\cdot) -\vartheta\Vert \partial_x^2 y(t) \Vert^2 +\vartheta\int_0^{\infty}\beta^{\prime} (s)\eta^2 \,ds-\nu_1 \int_{0}^{1} y^2(x,t) \partial_x y(x,t) dx,
\end{equation*} 
therefore, using the Dirichlet boundary conditions in (\ref{11})$_2$, we see that
\begin{equation*}
-\nu_1 \int_{0}^{1} y^2 (x,t) \partial_x y(x,t) dx=-\frac{\nu_1}{3} \left[y^3\right]_0^1 =0,
\end{equation*}    
thus the above two properties imply that
\begin{equation}\label{EstDec1}
E^{\prime}(t)\leq -\vartheta(\partial_{x}y)^2 (0,t) -\vartheta\Vert \partial_x^2 y(t) \Vert^2 +\vartheta\int_0^{\infty}\beta^{\prime} (s)\eta^2 \,ds.
\end{equation}  
By combining \eqref{EstDec1} with \eqref{wi1} and \eqref{wi2}, it follows that 
\begin{equation}\label{EstDec3}
E^{\prime}(t)\leq -\pi^4\vartheta\Vert y(\cdot,t) \Vert^2 = -2\pi^4\vartheta E(t)+\pi^4\vartheta\int_0^{\infty}\beta (s)\eta^2 \,ds. 
\end{equation}
It is clear that \eqref{EstDec3} is similar to \eqref{Energy1}, so it leads to \eqref{Energy3} with $\pi^4\vartheta$ instead of $\frac{\omega_0}{c_0}$. 
\vskip0,1truecm
If $\xi^{\prime} =0$, we see that \eqref{Energy30} holds with $\vartheta$ instead of $\kappa$, and then \eqref{Energy31} is valid. Consequently, we get the exponential decay estimate \eqref{decay1}.  
\vskip0,1truecm
If $\xi^{\prime} \ne 0$, and according to \eqref{EstDec1}, we observe that \eqref{Energy4} is valid with $\vartheta$ instead of $\kappa$. Therefore, the same computations show that \eqref{Energy50}, \eqref{Energy6}, \eqref{BoundBound} and \eqref{Energy7} are still valid. Consequently, the proof of \eqref{decay2} can be ended as in the proof of Theorem \ref{theorem1} - Subcase 1.2.  
\end{proof}
\vskip0,1truecm
\begin{remark}
The reader has certainly noticed that well-posedness result of the KS problem is established under the condition $\nu_i > 0$, for $i=0,2,3$ (see \eqref{ai1}). Of course,  $\nu_0$ and $\nu_2$ are positive as they represent respectively the viscosity term coefficient and the anti-diffusion parameter. Notwithstanding, the requirement $\nu_3 > 0$ is used for sake of simplicity, and hence, can be relaxed. In fact, one can assume that $\nu_3$ is a non positive real constant and then appropriate modifications should be made. For instance, using the trace inequality
\[ (\partial_{x}y)^2 (0,\cdot) \le 2 \Vert \partial_x y \Vert^2 +\Vert \partial_x^2 y \Vert^2,\]
along with \eqref{wi2}, the estimate \eqref{diks} leads to
\begin{eqnarray}
\langle{\mathcal{K}}\Lambda,\Lambda \rangle_{\mathcal{H}}
&\le& \left[\left(\vert 1-\nu_0 \vert\dfrac{\epsilon}{2}-\nu_0 \nu_3\right)\left( \dfrac{2}{\pi^2}+1\right)+\dfrac{\nu_2}{\pi^2}-\nu_0 \right]  \Vert \partial_x^2 y \Vert^2\nonumber\\
&& +\frac{1}{2}\left[ 1-\dfrac{\alpha_0}{\epsilon}\vert 1-\nu_0\vert \right] \int_0^{\infty}\beta^{\prime}(s)\eta^2 ds, \label{diksn1}
\end{eqnarray}
for any $\epsilon >0$ and $\nu_3\leq 0$. Keeping the first and third conditions in \eqref{ai1} unchanged, one should require that the parameters $\alpha_0$ and $\nu_3$ obey the following weaker conditions than \eqref{f3} and the second one in \eqref{ai1}, respectively: 
$$
\nu_3 >\dfrac{1}{2+\pi^2}\left(\dfrac{\nu_2}{\nu_0} -\pi^2\right)
$$ 
and
$$
\alpha_0 <\dfrac{2}{(1-\nu_0 )^2}\left(\frac{\pi^2\nu_0-\nu_2}{2+\pi^2}+\nu_0\nu_3\right)\quad\hbox{if}\,\,\nu_0 \ne 1,
$$  
so that we can choose $\epsilon$ (instead of \eqref{eps}) as follows
\[
\alpha_0 \vert 1-\nu_0\vert < \epsilon < \dfrac{2}{\vert 1-\nu_0\vert}\left[ \dfrac{1}{2+\pi^2}\left(\pi^2\nu_0- \nu_2\right)+\nu_0 \nu_3 \right].\]
In this case, the dissipativity of the operator $\mathcal{K}$ follows from \eqref{diksn1} since we have
\[
\langle{\mathcal{K}}\Lambda,\Lambda \rangle_{\mathcal{H}}
\le -\vartheta \Vert \partial_x^2 y \Vert^2+\vartheta \int_0^{\infty}\beta^{\prime}(s)\eta^2 ds ,
\]
where instead of \eqref{vart}, the positive constant $\vartheta$ is
\[
  \vartheta=\min \left\{\frac{1}{2}\left[ 1-\dfrac{\alpha_0}{\epsilon}\vert 1-\nu_0\vert \right], \nu_0-\dfrac{\nu_2}{\pi^2}-\left(\vert 1-\nu_0 \vert\dfrac{\epsilon}{2}-\nu_0 \nu_3\right)\left( \dfrac{2}{\pi^2}+1\right) \right\}.
\]
Thereafter, running on much the same lines as previously done with of course a number of minor changes, we can obtain similar results to those in Theorem \ref{th1ks}, Theorem \ref{th2ks}, Theorem \ref{kssol} and Theorem \ref{theorem5511}.
\end{remark}

\section{Numerical analysis of \eqref{1} and \eqref{11}}\label{sec6}
\subsection{Generalized scheme proposal.}\label{sec6.1}
In this section, we will present a numerical scheme that solves both \eqref{1} and \eqref{11}. For $M\in\mathbb{N}$, we will discretize the interval $[0,1]$ using $M+1$ equally separated nodes. Let us define $x_k = k\Delta x,\: k =0,1,\dots,M$ and $\Delta x := \frac{1}{M}$. For the time variable, and for $n\in\mathbb{N}$, let $t_n := n\Delta t;\: \Delta t \in (0,1)$. With this, we will write $y_k^n \approx y(x_k,t_n)$; that is, $y_k^n$ will be our numerical solution at $x = x_k$ and $t = t_n$. Define $y_k^{n+\frac{1}{2}} := \frac{y_k^{n+1}+y_k^n}{2}$.  Due to computational limitations, we will consider a bounded domain $[0,s_f]$ for the $s$ variable, discretized using $L$ points $s_i = i\Delta s,\: i =0,1,\dots,L-1$, for $\Delta s < 1$ given.
\vskip0,1truecm
Because we are dealing with both \eqref{1} and \eqref{11}, let us focus our attention on the following PDE:
\begin{equation}\label{pde_gen}
\partial_t y(x,t) + \alpha_1 \partial_x y(x,t) + \alpha_2 \partial_x^2 y(x,t) + \alpha_3 \partial_x^3 y(x,t) + \alpha_4 \partial_x^4 y(x,t)+ \alpha_5 y(x,t)\partial_x y(x,t) = 0
\end{equation}
where $\alpha_i \in \mathbb{R}, i=1,2,3,4$. Thus, when $\alpha_4 = 0$ we recover \eqref{1}, while \eqref{11} is obtained when $\alpha_1 = \alpha_3 = 0$.
\vskip0,1truecm
We will approximate \eqref{pde_gen} using a finite differences approach. To this end, let us define the vector space
\[\hat{X}_M = \{u=(u_0\: u_1\: \dots \: u_M)^T\in \mathbb{R}^{M+1}:u_0 = u_M = 0 \},\]
and the vector subspace
\[ X_M := \{v \in \mathbb{R}^{M-3} : u = (u_0\: \dots \: u_{M})^T \in \hat{X}_M, v = (u_2\:u_3 \: u_4 \: \dots \: u_{M-2})^T \}.\]
Let $y^n \in X_M$. The derivatives will be approximated using
\begin{align}
  \partial_t y(x_k,t_n) &\approx \dfrac{y^{n+1}_{k} - y^{n}_{k}}{\Delta t}; \nonumber \\
  \partial_x y(x_k,t_n) &\approx D_xy_k^n := \dfrac{y^n_{k+1} - y^n_{k-1}}{2\Delta x}; \label{dif_fin} \\
  \partial^2_x y(x_k,t_n) &\approx D^2_xy_k^n := \dfrac{y^n_{k+1} -2y_k^n +  y^n_{k-1}}{\Delta x^2}; \nonumber \\
  \partial^3_x y(x_k,t_n) &\approx D^3_xy_k^n := \dfrac{\frac{1}{2}y^n_{k+2} - y^n_{k+1} + y^n_{k-1} - \frac{1}{2}y^n_{k-2}}{\Delta x^3}; \nonumber \\
  \partial^4_x y(x_k,t_n) &\approx D^4_xy_k^n := \dfrac{y^n_{k+2} - 4y^n_{k+1} + 6y_k^n -4y^n_{k-1} + y^n_{k-2}}{\Delta x^4}. \nonumber
\end{align}
All these are second-order approximations of their respective derivatives. They also induce the definition of the following matrix operators in $\mathbb{R}^{(M-3)\times (M-3)}$ over $y^n$:
\begin{align*}
  \boldsymbol{D}y^n := \dfrac{1}{\Delta x}\begin{pmatrix}0 & \frac{1}{2}& & & \\ -\frac{1}{2} & 0 & \frac{1}{2} & & \\ & \ddots & \ddots & \ddots & \\ & & -\frac{1}{2} & 0 & \frac{1}{2}  \\  & & & -\frac{1}{2} & 0  \end{pmatrix},\:
  &\boldsymbol{D}_x^2y^n := \dfrac{1}{\Delta x^2}\begin{pmatrix}-2 & 1& & & \\ 1 & -2 & 1 & & \\ & \ddots & \ddots & \ddots & \\ & & 1 & -2 & 1  \\  & & & 1 & -2  \end{pmatrix} \\
  \boldsymbol{D}_x^3y^n := \dfrac{1}{\Delta x^3}\begin{pmatrix}0 & -1& \frac{1}{2} & & & & \\ 1 & 0 & -1 & \frac{1}{2} & & & \\ -\frac{1}{2} & 1 & 0 & -1 & \frac{1}{2} & & \\ & \ddots & \ddots & \ddots & \ddots & \ddots & \\ & & -\frac{1}{2} & 1 & 0 & -1 & \frac{1}{2}  \\ & & & -\frac{1}{2} & 1 & 0 & -1 \\ & & & & -\frac{1}{2} & 1 & 0  \end{pmatrix},\:
  &\boldsymbol{D}_x^4y^n := \dfrac{1}{\Delta x^4}\begin{pmatrix}6 & -4& 1 & & & & \\ -4 & 6 & -4 & 1 & & & \\ 1 & -4 & 6 & -4 & 1 & & \\ & \ddots & \ddots & \ddots & \ddots & \ddots & \\ & & 1 & -4 & 6 & -4 & 1  \\ & & & 1 & -4 & 6 & -4 \\ & & & & 1 & -4 & 6  \end{pmatrix}.
\end{align*}
Thus, our generalized Crank-Nicholson numerical scheme for \eqref{pde_gen} will be defined as follows: find $y^{n+1} \in X_M$ such that
\begin{equation}\label{esq_num}
  \dfrac{y^{n+1}-y^n}{\Delta t} + \alpha_1 \boldsymbol{D}_xy^{n+\frac{1}{2}} + \alpha_2 \boldsymbol{D}^2_xy^{n+\frac{1}{2}} + \alpha_3 \boldsymbol{D}^3_xy^{n+\frac{1}{2}} + \alpha_4 \boldsymbol{D}^4_xy^{n+\frac{1}{2}} + \alpha_5 y^{n+\frac{1}{2}}\boldsymbol{D}_xy^{n+\frac{1}{2}} = 0,\: n\in\mathbb{N},
\end{equation}
for $y^0\in X_M$ given. The boundary conditions will be considered in subsections \ref{sec6.2} and \ref{sec6.3}. In order to solve for $y^{n+1}$, we will have to solve the following problem for each timestep:
\begin{equation}\label{fixed_p}
    \begin{matrix}
  \left(\boldsymbol{I} + \alpha_1 \dfrac{\Delta t}{2} \boldsymbol{D}_x + \alpha_2 \dfrac{\Delta t}{2}\boldsymbol{D}^2_x + \alpha_3 \dfrac{\Delta t}{2}\boldsymbol{D}^3_x + \alpha_4 \dfrac{\Delta t}{2} \boldsymbol{D}^4_x\right) y^{n+1} \\
  =   \left(\boldsymbol{I} - \alpha_1 \dfrac{\Delta t}{2} \boldsymbol{D}_x - \alpha_2 \dfrac{\Delta t}{2}\boldsymbol{D}^2_x - \alpha_3 \dfrac{\Delta t}{2}\boldsymbol{D}^3_x - \alpha_4 \dfrac{\Delta t}{2} \boldsymbol{D}^4_x\right) y^{n} - \alpha_5 y^{n+\frac{1}{2}}\boldsymbol{D}_x y^{n+\frac{1}{2}} - \boldsymbol{f}^{n} - \boldsymbol{f}^{n+1}
  \end{matrix},
\end{equation}
where $\boldsymbol{I} \in \mathbb{R}^{(M+1)\times (M+1)}$ is the identity matrix. Inturn, the vectors $\boldsymbol{f^n}\in\mathbb{R}^{M-3}$, containing the boundary terms, will be properly defined later. This is a nonlinear problem which will be solved using a Picard fixed point iteration. This means that we need to solve a pentadiagonal system of equations many times per timestep. Because the structure of the coefficient matrix is the same during the whole simulation, an LU decomposition is computed only once using the \texttt{LAPACK}\footnote{https://www.netlib.org/lapack/} package for \texttt{FORTRAN 90}, and used to solve for the rest of the calculations. \\

 We will explain why the matrix operators do not act directly over $y_1^n$ and $y_{M-1}^n$. As this will be considered for both the KdVB and KS equations, let us pay our attention to the approximation of the fourth derivative at $x = x_2$, $x=x_3$, $x = x_{M-3}$, and $x = x_{M-2}$:
\[D_x^4 y_2^{n+1} = \dfrac{y^{n+1}_0-4y_1^{n+1}+6y^{n+1}_2-4y^{n+1}_3+y^{n+1}_4}{\Delta x^4}; \]
\[D_x^4 y_3^{n+1} = \dfrac{y^{n+1}_1-4y_2^{n+1}+6y^{n+1}_3-4y^{n+1}_4+y^{n+1}_5}{\Delta x^4}; \]
\[D_x^4 y_{M-3}^{n+1} = \dfrac{y^{n+1}_{M-5}-4y^{n+1}_{M-4}+6y^{n+1}_{M-3}-4y^{n+1}_{M-2}+y^{n+1}_{M-1}}{\Delta x^4}; \]
\[D_x^4 y_{M-2}^{n+1} = \dfrac{y^{n+1}_{M-4}-4y^{n+1}_{M-3}+6y^{n+1}_{M-2}-4y^{n+1}_{M-1}+y^{n+1}_M}{\Delta x^4}. \]
Since $y_0^{n+1} = y_M^{n+1} = 0$, and $y_1^{n+1}$ with $y_{M-1}^{n+1}$ are both known, we will define the matrix operator $\boldsymbol{\hat{D}}_x^4$ as follows
\[\boldsymbol{\hat{D}}_x^4 y^{n+1} = \dfrac{1}{\Delta x^4}\begin{pmatrix}6 & -4& 1 & & & & \\ -4 & 6 & -4 & 1 & & & \\ 1 & -4 & 6 & -4 & 1 & & \\ & \ddots & \ddots & \ddots & \ddots & \ddots & \\ & & 1 & -4 & 6 & -4 & 1  \\ & & & 1 & -4 & 6 & -4 \\ & & & & 1 & -4 & 6  \end{pmatrix} \begin{pmatrix} y_2^{n+1} \\ y_3^{n+1} \\ \vdots \\ \vdots \\ y_{M-3}^{n+1} \\ y_{M-2}^{n+1} \end{pmatrix} + \dfrac{1}{\Delta x^4} \begin{pmatrix}- 4y_1^{n+1} \\ y_1^{n+1} \\ 0 \\ \vdots \\ 0 \\ y_{M-1}^{n+1} \\ -4y_{M-1}^{n+1} \end{pmatrix};
\]
this is,
\[
\boldsymbol{\hat{D}}_x^4 y^{n+1} = \boldsymbol{D}_x^4 y^{n+1} + \dfrac{1}{\Delta x^4} \begin{pmatrix}- 4y_1^{n+1} \\ y_1^{n+1} \\ 0 \\ \vdots \\ 0 \\ y_{M-1}^{n+1} \\ -4y_{M-1}^{n+1} \end{pmatrix}.
\]
In a similar fashion, we will re-define the matrix operators for the other derivatives:
\[\boldsymbol{\hat{D}}_x y^{n+1} = \dfrac{1}{\Delta x} \begin{pmatrix}0 & \frac{1}{2}& & & \\ -\frac{1}{2} & 0 & \frac{1}{2} & & \\ & \ddots & \ddots & \ddots & \\ & & -\frac{1}{2} & 0 & \frac{1}{2}  \\  & & & -\frac{1}{2} & 0  \end{pmatrix}     \begin{pmatrix} y_2^{n+1}  \\ \vdots \\ \vdots  \\ y_{M-2}^{n+1} \end{pmatrix} + \dfrac{1}{\Delta x} \begin{pmatrix} -\frac{1}{2}y_1^{n+1} \\  0 \\ \vdots \\ 0  \\ \frac{1}{2}y_{M-1}^{n+1} \end{pmatrix},
\]
\[\boldsymbol{\hat{D}}_x^2 y^{n+1} = \dfrac{1}{\Delta x^2} \begin{pmatrix}-2 & 1& & & \\ 1 & -2 & 1 & & \\ & \ddots & \ddots & \ddots & \\ & & 1 & -2 & 1  \\  & & & 1 & -2  \end{pmatrix}   \begin{pmatrix} y_2^{n+1} \\  \vdots \\ \vdots  \\ y_{M-2}^{n+1} \end{pmatrix} + \dfrac{1}{\Delta x^2} \begin{pmatrix} y_1^{n+1} \\  0 \\ \vdots \\ 0  \\ y_{M-1}^{n+1} \end{pmatrix},
\]
\[\boldsymbol{\hat{D}}_x^3 y^{n+1} = \dfrac{1}{\Delta x^3} \begin{pmatrix}0 & -1& \frac{1}{2} & & & & \\ 1 & 0 & -1 & \frac{1}{2} & & & \\ -\frac{1}{2} & 1 & 0 & -1 & \frac{1}{2} & & \\ & \ddots & \ddots & \ddots & \ddots & \ddots & \\ & & -\frac{1}{2} & 1 & 0 & -1 & \frac{1}{2}  \\ & & & -\frac{1}{2} & 1 & 0 & -1 \\ & & & & -\frac{1}{2} & 1 & 0  \end{pmatrix}    \begin{pmatrix} y_2^{n+1} \\ y_3^{n+1} \\ \vdots \\ \vdots \\ y_{M-3}^{n+1} \\ y_{M-2}^{n+1} \end{pmatrix} + \dfrac{1}{\Delta x^3} \begin{pmatrix} y_1^{n+1} \\ -\frac{1}{2}y_1^{n+1} \\ 0 \\ \vdots \\ 0 \\ \frac{1}{2}y_{M-1}^{n+1} \\ -y_{M-1}^{n+1} \end{pmatrix}.
\]
This motivates the definition of the vector
\[
\boldsymbol{f}^{n+1} :=  \dfrac{1}{\Delta x} \begin{pmatrix} -\frac{1}{2}y_1^{n+1} \\  0 \\ \vdots \\ 0  \\ \frac{1}{2}y_{M-1}^{n+1} \end{pmatrix} + \dfrac{1}{\Delta x^2} \begin{pmatrix} y_1^{n+1} \\  0 \\ \vdots \\ 0  \\ y_{M-1}^{n+1} \end{pmatrix} +  \dfrac{1}{\Delta x^3} \begin{pmatrix} y_1^{n+1} \\ -\frac{1}{2}y_1^{n+1} \\ 0 \\ \vdots \\ 0 \\ \frac{1}{2}y_{M-1}^{n+1} \\ -y_{M-1}^{n+1} \end{pmatrix} + \dfrac{1}{\Delta x^4} \begin{pmatrix}- 4y_1^{n+1} \\ y_1^{n+1} \\ 0 \\ \vdots \\ 0 \\ y_{M-1}^{n+1} \\ -4y_{M-1}^{n+1} \end{pmatrix},
\]
which in turns leads to \eqref{fixed_p}. \\
\vskip0,1truecm
As our focus will be on the energy decay, we will define the following discrete analogue:
\[E(y^n) := \dfrac{1}{2}\sum\limits_{i=1}^{M-1}|y_i^n|^2\Delta x + \sum\limits_{i=0}^{L}\beta(s_i)\eta^2(t_n,s_i) \Delta s. \]

\subsection{Boundary conditions for \eqref{1}}\label{sec6.2}
While the boundary conditions for $y(x,t)$ at $x=0$ and $x=1$ are already imposed in the definition of the vector space $X_M$, we still need to consider the conditions for $\partial_x y(x,t)$. This means we need to impose conditions for $y_1^{n+1}$ and $y_{M-1}^{n+1}$. Let us recall the memory term associated to \eqref{1}
\[ \partial_x y(1,t)=\omega_4 \partial_x y(0,t)+\displaystyle \int_0^{\infty}\alpha (s) \, \partial_x y(0,t-s) ds. \]
Discretizing the terms outside the integral, we get
\[D_x y^{n+1}_M=\omega_4 D_x y^{n+1}_0+\displaystyle \int_0^{\infty}\alpha (s) \, \partial_x y(0,t-s) ds, \]
and after considering \eqref{dif_fin},
\[\dfrac{y_{M+1}^{n+1} - y_{M-1}^{n+1}}{2\Delta x}=\omega_4 \dfrac{y_{1}^{n+1} - y_{-1}^{n+1}}{2\Delta x}+\displaystyle \int_0^{\infty}\alpha (s) \, \partial_x y(0,t-s) ds. \]
Observe that extra nodes at $x = (M+1)\Delta x$ and $x = -\Delta x$ appear; instead, we will assume that $y_{M+1}^{n+1} = y_{-1}^{n+1} = 0$ due to our already known boundary conditions. Thus, the previous expression turns into
\begin{equation}\label{front1}
  0= \dfrac{y^{n+1}_{M-1}}{2\Delta x} + \omega_4 \dfrac{y_1^{n+1}}{2\Delta x}+\displaystyle \int_0^{\infty}\alpha (s) \, \partial_x y(0,t-s) ds.
\end{equation}
Therefore, we will only need to compute $y_1^{n+1}$ or $y_{M-1}^{n+1}$. Let us turn our attention to $y_1^{n+1}$, and let us recall \eqref{etatheta}
\[
\displaystyle\int_0^{\infty}\,\beta (s) \eta (t,s)\,ds=\displaystyle\int_0^{\infty}\,\alpha (s) \partial_x y (0,t-s) \,ds.
\]
Evaluating the latter at $t=t_{n+1}$ and separating the integral in the right hand side, we have
\[
\displaystyle\int_0^{\infty}\,\beta (s) \eta (t_{n+1},s)\,ds=\displaystyle\int_0^{t_{n+1}}\,\alpha (s) \partial_x y (0,t_{n+1}-s) \,ds + \int_{t_{n+1}}^{\infty}\,\alpha (s) \partial_x y (0,t_{n+1}-s) \,ds.
\]
Since we know the value of $\partial_x y$ only at discrete values of $t$, and the function $\alpha$ is previously known in its exact form, we will approximate the first integral of the right side of the last identity as follows
\[ \displaystyle\int_0^{t_{n+1}}\,\alpha (s) \partial_x y (0,t_{n+1}-s) \,ds \approx \displaystyle\sum_{i=0}^{n+1}\,\alpha (s_i) \partial_x y (0,t_{n+1}-s_i) \,\Delta t, \]
where, in {\it this integral}, $s_i = i\Delta t$. After considering this, and pulling out the first term in the finite sum, we can rewrite \eqref{etatheta} as
\[\displaystyle\int_0^{\infty}\,\beta (s) \eta (t_{n+1},s)\,ds = \displaystyle \alpha(0)\partial_x y (0,t_{n+1})\Delta t+ \sum_{i=1}^{n+1}\,\alpha (s_i) \partial_x y (0,t_{n+1}-s_i) \,\Delta t +  \int_{t_{n+1}}^{\infty}\,\alpha (s) \partial_x y (0,t_{n+1}-s) \,ds, \]
where the approximation sign was changed to an equality, because this is the expression we will manipulate to compute $y_1^{n+1}$. Recalling \eqref{dif_fin}, we get
\[\displaystyle\int_0^{\infty}\,\beta (s) \eta (t_{n+1},s)\,ds = \displaystyle \alpha(0)\dfrac{y_1^{n+1} - y_{-1}^{n+1}}{2\Delta x}\Delta t+ \sum_{i=1}^{n+1}\,\alpha (s_i) \partial_x y (0,t_{n+1}-s_i) \,\Delta t +  \int_{t_{n+1}}^{\infty}\,\alpha (s) \partial_x y (0,t_{n+1}-s) \,ds, \]
because $y_{-1}^{n+1} = 0$, and after re-arranging terms, we get an expression to compute $y_1^{n+1}$:
\begin{equation}\label{front2}
  \displaystyle  y_1^{n+1} = \displaystyle \dfrac{2\Delta x}{\alpha(0)\Delta t}\left( \int_0^{\infty}\,\beta (s) \eta (t_{n+1},s)\,ds - \sum_{i=1}^{n+1}\,\alpha (s_i) \partial_x y (0,t_{n+1}-s_i) \,\Delta t -  \int_{t_{n+1}}^{\infty}\,\alpha (s) \partial_x y (0,t_{n+1}-s) \,ds\right).
\end{equation}
Replacing in \eqref{front1}, we get $y_{M-1}^{n+1}$. The other integrals involved are computed using a Simpson's Rule.

%% The other integral involved can be computed using another discretization for the $s$ variable:
%% \[\displaystyle\int_0^{\infty}\,\beta (s) \eta (t_{n+1},s)\,ds \approx \sum_{j=0}^K\beta (s_j) \eta (t_{n+1},s_j)\Delta s,  \]
%% for $s_j = j\Delta s,\: \Delta s \in (0,1)$ and some $K \in \mathbb{N}$.

\subsection{Numerical experiments for the KdVB problem}\label{sec6.3}
\subsubsection{Case 1.}

\begin{figure}
  \includegraphics[scale=0.65]{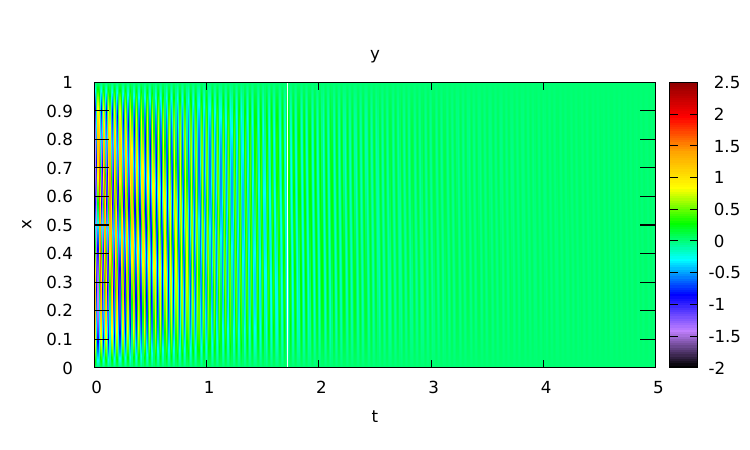}
  \includegraphics[scale=0.6]{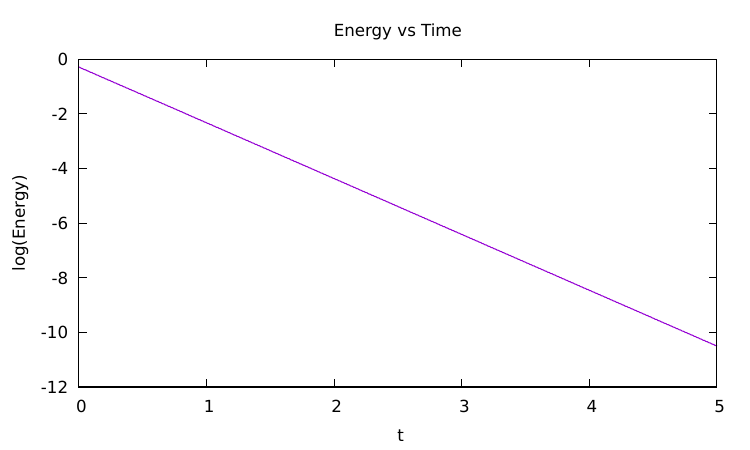}
  \caption{Results for Case 1.}
  \label{caso1}
\end{figure}

Let us present some results regarding \eqref{1}. For this case, we will use $w_0 = 0.01$, $w_1 = 1$, $w_2 = 2$, $w_3 = 6$, $w_4 = 0.1$; $t\in (0,5]$, $s\in [0,30]$; $\alpha(s) = d_2e^{-d_1 t}$, $d_1 = 2$ and $d_2 = 0.01$; $y_0(x) = 1-\cos(2\pi x)$; $y_1(t) = 0$, and thus, $\eta_0(s) = 0$; $\Delta s = 1.8311E-3$, $\Delta t = 4.8875E-3$, and $\Delta x = 1.22E-4$. Figure \ref{caso1} illustrates our numerical results. It is clear that the energy decay is exponential, as expected from the previous study.

\subsubsection{Case 2.}

Here we will consider a non-zero function $y_1$. $w_0 = 0.005$, $w_1 = 1$, $w_2 = 2$, $w_3 = 6$, $w_4 = -0.9$; $t\in (0,5]$, $s\in [0,30]$; $\alpha(s) = d_2e^{-d_1 t}$, $d_1 = 1$ and $d_2 = 0.01$; $y_0(x) = 1-\cos(2\pi x)$; $y_1(t) = \dfrac{\sin(t)}{10}$, and thus, $\eta_0(s) = \dfrac{1-\cos(s)}{10}$; $\Delta s = 1.831E-3$, $\Delta t = 4.887E-3$, and $\Delta x = 1.22E-4$. Figure \ref{caso1} illustrates our numerical results. We can see that the energy decay is exponential as well.

\begin{figure}
  \includegraphics[scale=0.65]{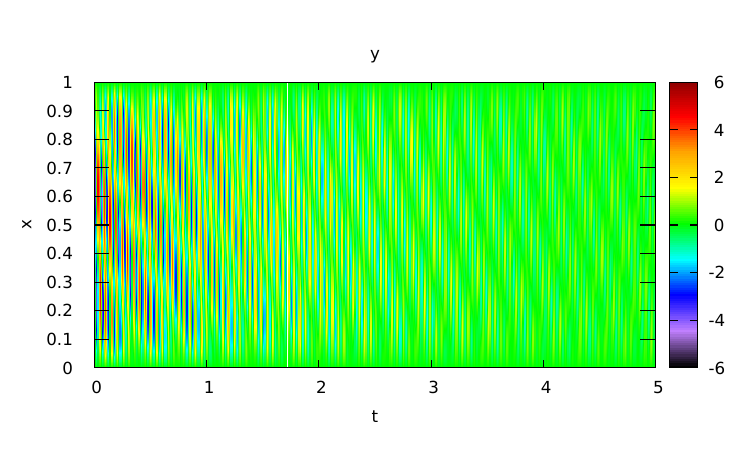}
  \includegraphics[scale=0.6]{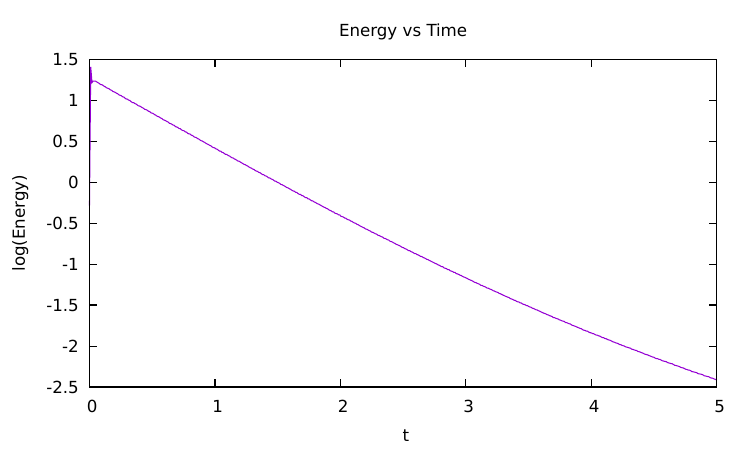}
  \caption{Results for Case 2.}
  \label{caso2}
\end{figure}

\subsection{Boundary conditions for \eqref{11}}\label{sec6.3}
Let us focus our attention now on problem \eqref{11}. The additional condition $\partial_x^2 y(1,t)=0$ can be translated to
\[\dfrac{y_{M-1}^{n+1}-2y_{M}^{n+1} + y_{M+1}^{n+1}}{\Delta x^2} = 0. \]
The fact that $y_{M}^{n+1} = 0$ motivates us to consider $y_{M+1}^{n+1} = 0$; thus, $y_{M-1}^{n+1} = 0$ as well. Regarding the memory condition, we have
\[ \dfrac{y_{-1}^{n+1} - 2y_{0}^{n+1} + y_{1}^{n+1}}{\Delta x^2} = \nu_3 \dfrac{y_{1}^{n+1} - y_{-1}^{n+1}}{2\Delta x} + \int\limits_0^\infty \alpha(s) \partial_x y(0,t-s)ds, \]
where $y_{0}^{n+1}$. Let us observe, however, that we can obtain conditions for both $y_{-1}^{n+1}$ and $y_{1}^{n+1}$. In fact, let us recall, from \eqref{front2}, that
\begin{align} 
  \displaystyle 
    \alpha(0)\dfrac{y_1^{n+1} - y_{-1}^{n+1}}{2\Delta x}\Delta t& = \displaystyle\int_0^{\infty}\,\beta (s) \eta (t_{n+1},s)\,ds \nonumber  \\ &\:- \sum_{i=1}^{n+1}\,\alpha (s_i) \partial_x y (0,t_{n+1}-s_i) \,\Delta t -  \int_{t_{n+1}}^{\infty}\,\alpha (s) \partial_x y (0,t_{n+1}-s) \,ds. \label{front3}
\end{align}
  From which we obtain
  \[y_{-1}^{n+1} = y_{1}^{n+1} - \dfrac{2\Delta x}{\alpha(0)\Delta t}\left(\int_0^{\infty}\,\beta (s) \eta (t_{n+1},s)\,ds - \sum_{i=1}^{n+1}\,\alpha (s_i) \partial_x y (0,t_{n+1}-s_i) \,\Delta t -  \int_{t_{n+1}}^{\infty}\,\alpha (s) \partial_x y (0,t_{n+1}-s) \,ds\right). \]
  Replacing in \eqref{front3},
  \small
  \begin{align*}
    y_{1}^{n+1}\left(1-\dfrac{\nu_3\Delta x-2}{\nu_3\Delta x + 2} \right) &= \dfrac{2\Delta x}{\alpha(0)\Delta t}\left(\int_0^{\infty}\,\beta (s) \eta (t_{n+1},s)\,ds - \sum_{i=1}^{n+1}\,\alpha (s_i) \partial_x y (0,t_{n+1}-s_i) \,\Delta t -  \int_{t_{n+1}}^{\infty}\,\alpha (s) \partial_x y (0,t_{n+1}-s) \,ds\right) \\
    &+ \dfrac{2\Delta x^2}{2+\nu_3 \Delta x}\int\limits_0^\infty \alpha(s) \partial_x y(0,t-s)ds,
  \end{align*}
  \normalsize
  which allows us to proceed as for the KdVB case.

  \subsection{Numerical experiments for the KS problem}\label{sec6.4}
  \subsubsection{Case 3.}
\begin{figure}
  \includegraphics[scale=0.65]{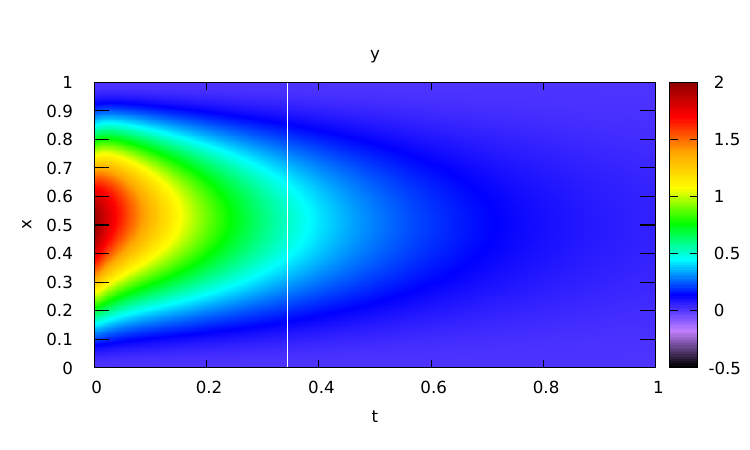}
  \includegraphics[scale=0.6]{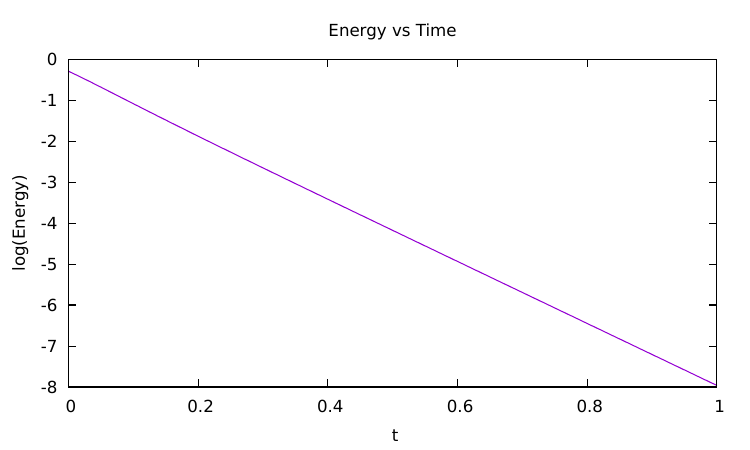}
  \caption{Results for Case 3.}
  \label{caso3}
\end{figure}

  As a first example, let us consider the KS equation with $\nu_0 = 0.01$, $\nu_1 = 1$, $\nu_2 = 0.1$, $\nu_3 = 0.1$; $t\in (0,1]$, $s\in[0,30]$; $\alpha(s) = d_2e^{-d_1 s}$ with $d_1 = 1$ and $d_2 = 0.1$; $y_1(t) = 0$; $y_0(x) = 1-\cos(2\pi x)$; $\Delta s =2.932E-2$, $\Delta t = 3.921E-3$ and $\Delta x = 9.737E-4$. Figure \ref{caso3} illustrates our numerical results.

\subsubsection{Case 4.} 
\begin{figure}
  \includegraphics[scale=0.65]{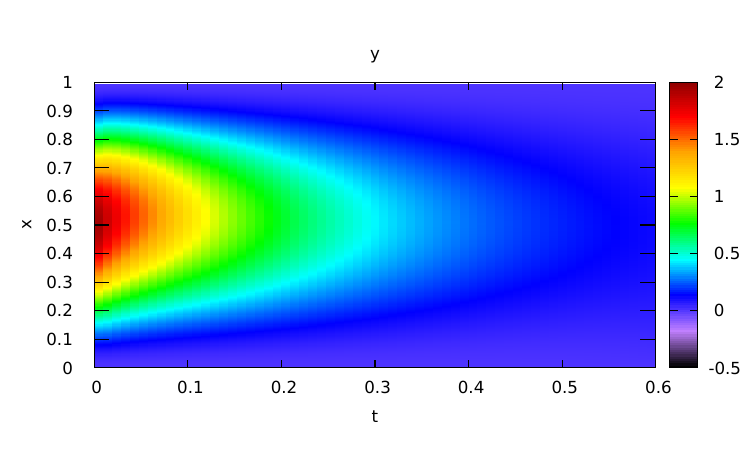}
  \includegraphics[scale=0.6]{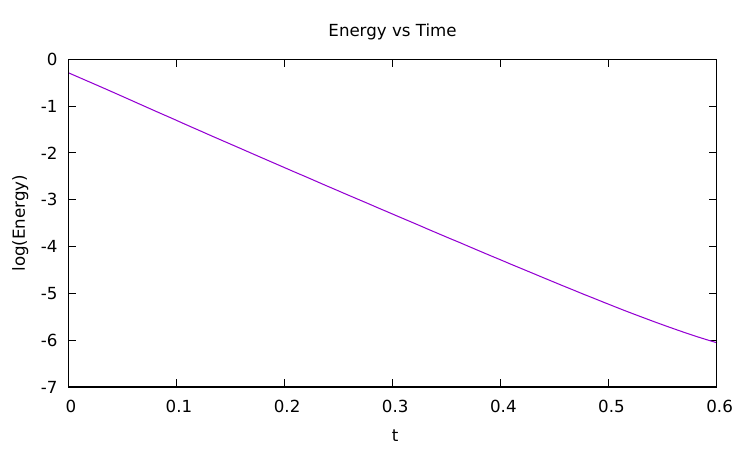}
  \caption{Results for Case 4.}
  \label{caso4}
\end{figure}

We will repeat Case 2 but using $y_1(t) = \dfrac{\sin(s)}{100}$; $t\in (0,0.6]$, $s\in [0,25]$; $\Delta s =2.443E-2$, $\Delta t = 9.523E-3$ and $\Delta x = 1.941E-3$. Results are in Figure \ref{caso4}.

\subsubsection{Case 5.} 
\begin{figure}
  \includegraphics[scale=0.8]{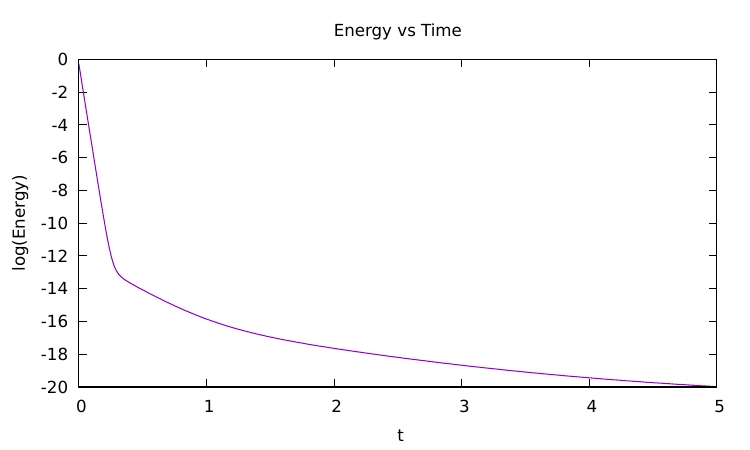}
  \caption{Energy decay in Case 5.}
  \label{caso5}
\end{figure}

As a final example, let us consider $\alpha(s) = d_2(1+s)^{-d_1}$, with $d_1 = 2$ and $d_2 = 0.01$. Regarding the other parameters, we will use $\nu_0 = 0.05$, $\nu_1 = 1$, $\nu_2 = 0.1$, $\nu_3 = 1$; $y_1(t) = 0$; $y_0(x) = 1-\cos(2\pi x)$; $t\in (0,5]$, $s\in [0,30]$; $\Delta s = 2.932E-2$, $\Delta t = 1.96E-2$ and $\Delta x = 9.737E-4$. Energy decay can be seen in Figure \ref{caso5}.

\section{Concluding discussion}\label{sec7}

In this paper, we provide an answer to the question posed in \cite{chengues2}. More precisely, we show that the well-posedness and stability properties for the KdVB equation are robust vis-\`a-vis a boundary memory-type control. Moreover, it shown that such a control contributes to the stability of the solutions. Obviously, our findings are obtained under some conditions on the physical parameters of the system, the memory kernel and the initial condition. This outcome is shown to be also true for the KdV equation, but more importantly, for a completely different type of problems related to the KS equation. Our results are ascertained by means of a numerical study.
\vskip0,1truecm
We aspire in a future work to investigate an interesting problem related to these KdVB and KS equations, not treated herein, is the well-posed and stability when the physical parameters of the equations vary either in terms of $x$ or $t$ or both. In the same spirit, the question what happens if the memory kernel is time-dependent is also paramount.

\section*{Acknowledgment} The numerical analysis was discussed during the visit of the third author to university of Lorraine, Metz, France, in March 2023, and the paper was finalized during the visit of the second author to University of Concepcion, Concepcion, Chile, in September 2023. The authors thank these two universities for their kind support and hospitality. 

\section*{Funding} MS was supported by  Fondecyt-ANID project 1220869, ANID (Chile) through project Centro de Modelamiento Matemático (BASAL projects ACE210010 and FB210005), ECOS-Sud project C20E03 (France - Chile),  INRIA Associated team ANACONDA, and Jean d'Alembert fellowship program, Université de Paris-Saclay.

\section*{Data availability statement} Data sharing is not applicable to the current paper as no data were generated or analysed during this study.

\section*{Conflict of Interest} The authors declare that they have no conflict of interest.

\end{document}